\documentclass[12pt,reqno]{amsart}
\usepackage[headings]{fullpage}
\usepackage{amssymb,amsmath,amscd,bbm,tikz}
\usepackage[all,cmtip]{xy}
\usepackage{url}
\usepackage[bookmarks=true,%
    colorlinks=true,%
    linkcolor=blue,%
    citecolor=blue,%
    filecolor=blue,%
    menucolor=blue,%
    urlcolor=blue,%
    breaklinks=true]{hyperref}
\usepackage{slashed}    
\usepackage{verbatim}
\usepackage{mathtools}
\usepackage[normalem]{ulem}
\graphicspath{{figures/}}

\newtheorem{theorem}{Theorem}[section]
\theoremstyle{definition}

\newtheorem{remark}[theorem]{Remark}
\newtheorem{corollary}[theorem]{Corollary}
\newtheorem{conjecture}[theorem]{Conjecture}

\def\BZ{\mathbbm Z}
\def\BQ{\mathbbm Q}

\def\BC{\mathbbm C}

\def\BE{\mathbbm E}

\def\calT{\mathcal T}

\def\ti{\widetilde}
\def\SL{\mathrm{SL}}

\def\ve{\varepsilon}

\def\be{\begin{equation}}
\def\ee{\end{equation}}

\def\CS{\mathrm{CS}}

\def\diag{\mathrm{diag}}

\begin{document}
\title[Twisted Neumann--Zagier matrices]{Twisted Neumann--Zagier matrices}
\author{Stavros Garoufalidis}
\address{
International Center for Mathematics, Department of Mathematics \\
Southern University of Science and Technology \\
Shenzhen, China \newline
{\tt \url{http://people.mpim-bonn.mpg.de/stavros}}}
\email{stavros@mpim-bonn.mpg.de}

\author{Seokbeom Yoon}
\address{Departament de Matem\`atiques \\
  Universitat Aut\`onoma de Barcelona \\
  Cerdanyola del Vall\`es, Spain \newline
{\tt \url{https://sites.google.com/view/seokbeom}}}
\email{sbyoon15@mat.uab.cat}



\date{30 August 2021}
\dedicatory{Dedicated to Walter Neumann and Don Zagier, with admiration.}

\begin{abstract}
  The Neumann--Zagier matrices of an ideal triangulation are integer matrices with
  symplectic properties whose entries encode the number of tetrahedra that wind around
  each edge of the triangulation. They can be used as input data for the construction
  of a number of quantum
  invariants that include the loop invariants, the 3D-index and state-integrals. 
  We define a twisted version of Neumann--Zagier matrices, describe their
  symplectic properties, and show how to compute them from the combinatorics of an ideal
  triangulation. As a sample application, we use them to define a twisted version of
  the 1-loop invariant (a topological invariant) which determines the 1-loop invariant
  of the cyclic covers of a hyperbolic knot complement, and conjecturally equals to the
  adjoint twisted Alexander polynomial.
\end{abstract}

\maketitle

{\footnotesize
\tableofcontents
}


\section{Introduction}
\label{sec.intro}

\subsection{Motivation}
\label{sub.motivation}

Ideal triangulations of 3-manifolds with torus boundary components were introduced by
Thurston~\cite{Thurston} as a convenient way to describe and effectively
compute~\cite{snappy} complete hyperbolic structures on 3-manifolds. To do so,
one assigns a complex number different from 0 or 1 to each tetrahedron and
a polynomial equation around each edge of the triangulation. These so-called gluing
equations have special shape that can be described by two matrices $\mathbf{A}$ and
$\mathbf{B}$ (with rows and columns indexed by the edges and by the tetrahedra,
respectively) whose entries describe the number of times (but \emph{not}
the order by which) tetrahedra wind around an edge.
One of the main discoveries of Neumann--Zagier is that the matrix
$(\mathbf{A} | \mathbf{B})$ becomes, after some minor modifications, the upper part
of a symplectic matrix with integer entries~\cite{NZ}. The symplectic property of the
NZ matrices and of the corresponding gluing equations define a linear symplectic
structure on a vector space whose quantization leads to a plethora of quantum
invariants that include the loop invariants of Dimofte and the first
author~\cite{DG1,DG2}, the 3D-index in both the original formulation of
Dimofte-Gaiotto-Gukov~\cite{DGG1,DGG2} as well as the state-integral formulation
of Kashaev and the first author~\cite{GK:meromorphic} and Kashaev--Luo--Vartanov
state-integral~\cite{KLV,AGK:KLV}. All of those invariants are defined using the
NZ matrices of a suitable ideal triangulation, and their topological invariance
follows by proving that they are unchanged under Pachner 2--3 moves.

Our original motivation was to study the behavior of the loop invariants of~\cite{DG1}
under cyclic covers. Since the latter are defined in terms of NZ data of an ideal
triangulation, we were led to study the behavior of the NZ matrices under cyclic
covers. By elementary topology, each tetrahedron of an ideal triangulation lifts to
$n$ tetrahedra in the $n$-fold cyclic cover, and lifting all the way to the
infinite cyclic cover leads to the notion of NZ matrices which we can call
\emph{equivariant} or \emph{twisted} (as is common in algebraic and geometric topology)
or \emph{$t$-deformed} (as is common in physics). We will use the term ``twisted'',
and keep in mind that the variable $t$ below encodes topological information of cyclic
covers.

Said differently, as the triangulation of the cyclic cover unfolds, so do its edges
and the tetrahedra that wind around them. This is the content of the twisted NZ
matrices. How can this elementary idea be nontrivial or interesting?

It turns out
that the twisted NZ matrices have twisted symplectic properties which come from topology
and 
using them one can give twisted versions of the above mentioned
invariants, i.e., of the loop invariants~\cite{DG1,DG2}, the 3D-index
in~\cite{DGG1,DGG2} and~\cite{GK:meromorphic} and KLV state-integral~\cite{KLV,AGK:KLV}.
A key property of such a twisted invariant is that it determines the corresponding
(untwisted) invariant of all cyclic covers. For example, the twisted 1-loop
invariant of a knot complement defined below satisfies this property.

Our goal is to define the twisted NZ matrices, describe their properties and show how
to compute them in terms of the methods developed by \texttt{SnapPy}~\cite{snappy}.
Having done so, we can use the twisted NZ matrices to define a twisted version of
the 1-loop invariant, prove its topological invariance and conjecture that it
equals to the adjoint twisted Alexander
polynomial. It is interesting to note that the twisted 1-loop invariant depends only
on the combinatorics of the NZ matrices of the infinite cyclic cover (which is
abelian information) whereas the adjoint twisted Alexander polynomial depends on the 
complete hyperbolic structure given as a representation of the fundamental group.

Further applications of twisted NZ matrices will be given in forthcoming work. 

\subsection{Torsion and its twisted version}
\label{sub.tt}

Before discussing twisting matters further, let us recall a key motivating example.
A basic invariant of a compact 3-manifold is the order of the torsion of its
first homology. The behavior of this invariant for all cyclic covers of a knot
complement is determined by a single Laurent polynomial, the Alexander polynomial.
Explicitly, we have
\be
\label{TDn}
\left|\mathrm{tor}(H_1(M^{(n)};\BZ))\right| = \prod_{\omega^n=1} \Delta_K(\omega) 
\ee
where $\Delta_K(t) \in \BZ[t^{\pm1}]$ is the  Alexander polynomial of a knot
$K \subset S^3$, $M^{(n)}$ is the $n$-fold cyclic cover of $M=S^3\setminus K$,
and the right-hand side
is the order of the torsion part of $H_1(M^{(n)};\BZ)$. This classical result
connecting the torsion of the first homology to the Alexander polynomial
(see e.g. \cite[(6.3), p.417]{Fox3}) is deeply rooted in the idea that the Alexander
polynomial is the torsion of the infinite cyclic cover of $M$ twisted by the
abelianization map $\alpha:\pi_1(M) \rightarrow H_1(M;\BZ)=\BZ$. In other words,
we have
\be
\boxed{
\label{basic}
\text{The twisted torsion determines the torsion of the cyclic covers} 
}
\ee

This idea has been extended in several directions.
Among them, one can define the torsion of a 3-manifold using interesting
representations of its fundamental group. For example, when $M$
is a cusped hyperbolic 3-manifold, one can define a torsion using a (lifted) geometric
$\SL_2(\BC)$-representation $\rho$ of its complement, or a symmetric power thereof.
There are some technical difficulties that one must overcome, stemming from the fact
that sometimes the corresponding chain complexes are not acyclic, hence the
torsion depends on a choice of peripheral curves, as well as normalization issues, since
the torsion is usually well-defined up to a sign. These issues have been addressed
in detail by~\cite{Porti:torsion,DubYam,Dunfield:twisted}. Among the several torsion
invariants,
we will be interested in the \emph{adjoint (Reidemeister) torsion} $\tau_\gamma(M)$
of a one-cusped hyperbolic 3-manifold (such as a knot complement) using the adjoint
representation $\mathrm{Ad}_\rho=\text{Sym}^2(\rho): \pi_1(M) \to \SL_3(\BC)$, where
$\gamma$ is a fixed peripheral curve. The adjoint torsion $\tau_\gamma(M) \in F^\times/\pm$
is a nonzero element of the trace field $F$ of $M$, well-defined up
to a sign~\cite{Porti:torsion}. Just as in the case of the Alexander polynomial,
there is a version of the \emph{adjoint twisted Alexander polynomial}
$\tau(M,\alpha,t) \in F[t^{\pm1}]/(\pm t^\BZ)$ (abbreviated by $\tau(M,t)$ when
$\alpha$ is clear) defined in~\cite{Wad94,Dunfield:twisted,DubYam} using an epimorphism
$\alpha: \pi_1(M) \to \BZ$. Here, the ambiguity is given by multiplication by an
element of $\pm t^\BZ:=\{\pm t^r \, | \, r \in \BZ\}$.

Two key properties of the adjoint twisted Alexander polynomial are
the behavior under finite cyclic covers~\cite{DubYam}
\be
\label{key1}
\tau(M^{(n)},t^n)=\prod_{\omega^n=1}\tau(M,\omega\, t)
\ee
and the relation with the adjoint torsion~\cite{Yamaguchi}, namely
\be
\label{key2}
\tau(M,1)=0, \qquad \left. \frac{d}{dt}\right|_{t=1}
\tau(M,t) = \tau_\lambda(M)
\ee
where $\lambda$ is the canonical longitude, i.e., the peripheral curve satisfying
$\alpha(\lambda)=0$. 

\subsection{Neumann--Zagier matrices and their twisted version}
\label{sub.tNZ}

\begin{figure}[htpb!]
\begingroup%
  \makeatletter%
  \providecommand\color[2][]{%
    \errmessage{(Inkscape) Color is used for the text in Inkscape, but the package 'color.sty' is not loaded}%
    \renewcommand\color[2][]{}%
  }%
  \providecommand\transparent[1]{%
    \errmessage{(Inkscape) Transparency is used (non-zero) for the text in Inkscape, but the package 'transparent.sty' is not loaded}%
    \renewcommand\transparent[1]{}%
  }%
  \providecommand\rotatebox[2]{#2}%
  \newcommand*\fsize{\dimexpr\f@size pt\relax}%
  \newcommand*\lineheight[1]{\fontsize{\fsize}{#1\fsize}\selectfont}%
  \ifx\svgwidth\undefined%
    \setlength{\unitlength}{114.47727954bp}%
    \ifx\svgscale\undefined%
      \relax%
    \else%
      \setlength{\unitlength}{\unitlength * \real{\svgscale}}%
    \fi%
  \else%
    \setlength{\unitlength}{\svgwidth}%
  \fi%
  \global\let\svgwidth\undefined%
  \global\let\svgscale\undefined%
  \makeatother%
  \begin{picture}(1,0.84370128)%
    \lineheight{1}%
    \setlength\tabcolsep{0pt}%
    \put(0,0){\includegraphics[width=\unitlength,page=1]{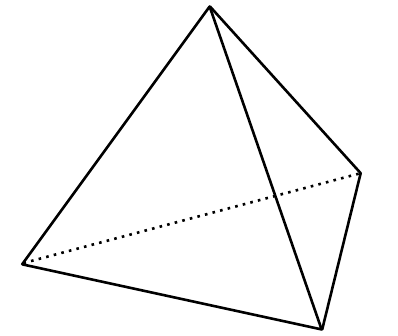}}%
    \put(0.22612007,0.52498993){\makebox(0,0)[lt]{\lineheight{1.25}\smash{\begin{tabular}[t]{l}$z$\end{tabular}}}}%
    \put(0.88209508,0.18215681){\makebox(0,0)[lt]{\lineheight{1.25}\smash{\begin{tabular}[t]{l}$z$\end{tabular}}}}%
    \put(0.72426566,0.63173901){\makebox(0,0)[lt]{\lineheight{1.25}\smash{\begin{tabular}[t]{l}$z''$\end{tabular}}}}%
    \put(0,0){\includegraphics[width=\unitlength,page=2]{tetrahedron.pdf}}%
    \put(0.60780402,0.46139955){\makebox(0,0)[lt]{\lineheight{1.25}\smash{\begin{tabular}[t]{l}$z'$\end{tabular}}}}%
    \put(0,0){\includegraphics[width=\unitlength,page=3]{tetrahedron.pdf}}%
    \put(0.45336339,0.2695276){\makebox(0,0)[lt]{\lineheight{1.25}\smash{\begin{tabular}[t]{l}$z'$\end{tabular}}}}%
    \put(0,0){\includegraphics[width=\unitlength,page=4]{tetrahedron.pdf}}%
    \put(0.42299074,0.07344528){\color[rgb]{0,0,0}\makebox(0,0)[lt]{\lineheight{1.25}\smash{\begin{tabular}[t]{l}$z''$\end{tabular}}}}%
  \end{picture}%
\endgroup%

\caption{An ideal tetrahedron.}
\label{fig.tetrahedron}
\end{figure}

Having already discussed the torsion and twisted torsion, we now recall some
basic facts from ideal triangulations of 3-manifolds and their gluing equations,
following~\cite{Thurston,NZ}. 
Let $M$ be an oriented hyperbolic 3-manifold with a torus boundary component
(often called a one-cusped hyperbolic 3-manifold) equipped with an ideal triangulation
$\calT$. An Euler characteristic argument shows that the number $N$ of tetrahedra in
$\calT$ is equal to the number of edges. We order the edges $e_i$ and the tetrahedra 
$\Delta_j$ of $\calT$ for $1 \leq i,j \leq N$. A quad type
of a tetrahedron is a pair of opposite edges; hence each tetrahedron has three
quad types. We fix an
orientation and a quad type of each tetrahedron $\Delta_j$ so that each edge of
$\Delta_j$ admits a shape parameter among 
\[z_j, \ z'_j=\frac{1}{1-z_j},\
  \textrm{or } z''_j=1- \frac{1}{z_j} \in \BC \setminus \{0,1\}\]
with opposite edges having same parameters (see Figure \ref{fig.tetrahedron}). 
We denote by $\square$ the three pairs of opposite edges of a tetrahedron,
so that the edges of $\square$ are assigned the edge parameter $z^\square$. 

A complete hyperbolic structure of $M$ can be described by a special solution
of the gluing equations. The latter describe the holonomy of the hyperbolic structure
around each edge $e_i$ of $\calT$, and have the (logarithmic) form
\be
\label{eqn:gluing}
\sum_{j=1}^N \left(\mathbf{G}_{ij} \log z_j +\mathbf{G}'_{ij} \log z'_j
  +\mathbf{G}''_{ij} \log z''_j\right) = 2\pi \sqrt{-1} .
\ee 
Here $\mathbf{G}$, $\mathbf{G}'$ and $\mathbf{G}''$ are the gluing equation matrices 
whose rows and columns are indexed by the edges and by the tetrahedra of $\calT$,
respectively, such that the $(i,j)$-entry of  $\mathbf{G}^\square$ is the number
of edges of $\Delta_j$ with parameter $z^\square_j$ is incident to the edge $e_i$
in $\calT$. It will be convenient to introduce $\zeta$-variables 
\be
\label{eqn.zeta}
\zeta = \frac{d \log z}{dz}=\frac{1}{z}, \quad
\zeta' = \frac{d \log z'}{dz}=\frac{1}{1-z},\quad
\zeta'' = \frac{d \log z''}{dz}=\frac{1}{z(z-1)} 
\ee
following \cite{siejakowski}. Note that the three shape parameters in each tetrahedron
satisfy the relation $z z' z''=-1$, and this implies the linear relation
$\zeta+\zeta'+\zeta''=0$.
Thus one can eliminate the variable $\zeta'$ from any expression given in terms of
$\zeta, \zeta'$, and $\zeta''$, arriving at the \emph{Neumann--Zagier matrices}
\be
\label{AB}
\mathbf{A}:=\mathbf{G}-\mathbf{G}', \qquad
\mathbf{B}:=\mathbf{G}''-\mathbf{G}'  .
\ee
Note that
\be
\label{3GAB}
\mathbf{G}\, \mathrm{diag}(\zeta)
+\mathbf{G}' \, \mathrm{diag}(\zeta')
+ \mathbf{G}'' \, \mathrm{diag}(\zeta'')
= \mathbf{A}  \, \mathrm{diag}(\zeta) + \mathbf{B} \, \mathrm{diag}(\zeta'') 
\ee
where $\mathrm{diag}(\zeta^\square)$ denotes the diagonal
matrix with diagonal entries $\zeta^\square_1,\ldots,\zeta^\square_N$.

The NZ matrices have remarkable properties discovered in~\cite{NZ}
which is a starting point for the quantization of the gluing equations and for
passing from hyperbolic geometry to quantum topology.
The symplectic property of the NZ matrices implies that 
\be 
\label{eqn.sympl}
\mathbf{A} \mathbf{B}^T = \mathbf{B} \mathbf{A}^T.
\ee

There are two enhancements of the gluing equations and of their corresponding
matrices. The first
is obtained by looking at the cusp, i.e., the peripheral (also called boundary)
torus of the 3-manifold.
By a  \emph{peripheral curve} we mean an oriented, homotopically non-trivial,
simple closed curve in the peripheral torus of $M$. A peripheral curve
$\gamma$ gives a triple $(\mathbf{C}_\gamma,\mathbf{C}'_\gamma,\mathbf{C}''_\gamma)$
of row vectors in $\BZ^N$ that describe the completeness equation as
\be
\sum_{j=1}^N \left(\mathbf{C}_{\gamma j}\log z_j +\mathbf{C}'_{\gamma j} \log z'_j
  +\mathbf{C}''_{\gamma j} \log z''_j\right)  = 0.
\ee
Fixing a peripheral curve $\gamma$, let $\widehat{\mathbf{G}}^\square$ denote
the matrix obtained from $\mathbf{G}^\square$  by replacing the last row by
$\mathbf{C}^\square_\gamma$, and likewise for $\widehat{\mathbf{A}}$ and
$\widehat{\mathbf{B}}$.

The second enhancement is a \emph{combinatorial flattening} of $\calT$, that
is a triple $(f,f',f'')$ of column vectors in $\BZ^N$ such that
\begin{align}
\mathbf{G} f+\mathbf{G}' f' + \mathbf{G}'' f'' & =(2,\ldots,2)^T,\\
f+ f' +  f'' & =(1,\ldots,1)^T,\\
\mathbf{C}_\gamma f+\mathbf{C}'_\gamma f' + \mathbf{C}''_\gamma f'' &=0	 
\end{align}
for any peripheral curve $\gamma$. This term was introduced in~\cite[Sec.4.4]{DG1}
as a necessary ingredient to define the 1-loop invariant (and there, it 
was called a combinatorial flattening compatible with any peripheral curve).
Every ideal triangulation has combinatorial flattenings~\cite[Thm.4.5]{Neumann04}.

This concludes our discussion of the NZ matrices. We now introduce a twisted version of
the NZ matrices of an ideal triangulation $\calT$ of a 3-manifold
$M$ as above. We will fix an epimorphism $\alpha : \pi_1(M) \rightarrow \BZ$ and a
peripheral curve $\mu$ satisfying $\alpha(\mu)=1$. For instance, if $M$ is a
knot complement in $S^3$, then $\alpha$ is the abelianization map and $\mu$ is a
meridian of the knot. Let $\ti M$ and $M^{(n)}$ denote the cyclic covers of $M$
corresponding to $\alpha^{-1}(0)$ and $\alpha^{-1}(n \BZ)$, respectively.
The ideal triangulation $\ti \calT$ of the infinite cyclic cover $\ti M$ induced
from $\calT$ is equipped with an action of the deck transformation group $\BZ$
(generated by the $\mu$-action), and induces an ideal triangulation $\calT^{(n)}$
of  the $n$-fold cyclic cover $M^{(n)}$.

We choose lifts $\ti e_i$ and $\ti \Delta_j$
of $e_i$ and $\Delta_j$ to $\ti \calT$, respectively, and for
$k \in \BZ$ let $\mathbf{G}_k^\square$ be the
$N \times N$ matrix 
whose $(i,j)$-entry is the number of edges of $\mu^k \cdot \ti \Delta_j$ with
parameter $z_j^\square$ is incident to the edge $\ti e_i$ in $\ti \calT$.
Note that $\mathbf{G}_k^\square$ is a zero matrix for all but finitely many $k$, as
there are finitely many tetrahedra around $\ti e_i$.
We define the \emph{twisted gluing equation matrices} $\mathbf{G}^\square(t)$ of $\calT$
by
\be
\label{Gt}
\mathbf{G}^\square(t):= \sum_{k \in \BZ} \mathbf{G}^\square_k\, t^k 
\ee
and the \emph{twisted Neumann--Zagier matrices} by
\be
\label{ABt}
\mathbf{A}(t):=\mathbf{G}(t)-\mathbf{G}'(t), \qquad
\mathbf{B}(t):=\mathbf{G}''(t)-\mathbf{G}'(t)  .
\ee
(More generally, if $\mathbf{X}(t)$ is a matrix with entries in $\BZ[t^{\pm1}]$,
we denote by $\mathbf{X}_k$ the coefficient of $t^k$ in $\mathbf{X}(t)$.) 
Note that the rows and columns of these matrices are indexed by the edges and
tetrahedra of $\calT$, respectively, and that their entries are in the ring
$\BZ[t^{\pm1}]$. Since the above matrices are well-defined after fixing lifts
of each edge of $\calT$, a different choice of lifts changes $\mathbf{G}^\square(t)$, 
$\mathbf{A}(t)$ and $\mathbf{B}(t)$ by multiplication from the left by the same
diagonal matrix $\diag(t^{c_1},t^{c_2},\dots,t^{c_N})$ for integers $c_1,c_2,\dots,c_N$.
This ambiguity propagates to any invariant constructed using these matrices.

A first property of the twisted NZ matrices $(\mathbf{A}(t) \,| \, \mathbf{B}(t))$
is that they determine the NZ matrices $(\mathbf{A}^{(n)} | \, \mathbf{B}^{(n)})$
of the cyclic covers $\calT^{(n)}$ of $\calT$, and moreover they do so via block
circulant matrices, i.e., square matrices such that each row is obtained by its
predecessor by a cyclic shift. We refer to~\cite{Davis} and~\cite{circulant} for
details on the properties of block circulant matrices.

\begin{theorem}
\label{thm.NZcyclic}
We have
\be
\label{An}
\mathbf{X}^{(n)} = 
\begin{pmatrix*}[c]
  \sum_{r \equiv 0} \mathbf{X}_r & \sum_{r \equiv 1} \mathbf{X}_r
  & \cdots & \sum_{r \equiv n-1} \mathbf{X}_r\\
  \sum_{r \equiv n-1} \mathbf{X}_r & \sum_{r \equiv 0} \mathbf{X}_r
  & \cdots & \sum_{r \equiv n-2} \mathbf{X}_r\\
  \vdots & \vdots & \ddots & \vdots \\
  \sum_{r \equiv 1} \mathbf{X}^\square_r & \sum_{r \equiv 2} \mathbf{X}_r
  & \cdots & \sum_{r \equiv 0} \mathbf{X}_r 
\end{pmatrix*}
\ee
for $\mathbf{X}=\mathbf{G}^\square$, $\mathbf{A}$ or $\mathbf{B}$.
\end{theorem}  

The twisted NZ matrices satisfy a twisted version of the symplectic
property of Equation~\eqref{eqn.sympl}.

\begin{theorem}
\label{thm.sympl}
We have
\be 
\label{eqn.abba}
\mathbf{A}(t) \mathbf{B}(1/t)^T = \mathbf{B}(t) \mathbf{A}(1/t)^T. 
\ee
\end{theorem}
In other words, $\mathbf{A}(t) \mathbf{B}(1/t)^T$ is symmetric under tranposition
followed by $t \mapsto 1/t$.
In particular, $\mathbf{A}(\omega) \mathbf{B}(\omega)^\ast$ is Hermitian for
$\omega \in \BC$ with $|\omega|=1$
and so is $\mathbf{B}(\omega)^{-1}\mathbf{A}(\omega)$ if $\mathbf{B}(\omega)$
is non-singular.
Here  $X^\ast$ is
the complex conjugate of the transpose of a matrix $X$.
 Note that given any ideal triangulation, we can choose quads
so that $\mathbf{B}$ (hence also $\mathbf{B}(t)$, since $\mathbf{B}(1)=\mathbf{B}$)
is invertible; see~\cite[App.A]{DG1}.

We now come to a conjectural property of the twisted NZ
matrices, which we have checked in numerous examples.

\begin{conjecture}
\label{conj.ABt}
The Laurent polynomials $\det(\mathbf{A}(t))$ and $\det(\mathbf{B}(t))$ are palindromic,
i.e. satisfy $p(t)=\ve t^r p(1/t)$ for some $\ve = \pm1$ and an integer $r$. 
\end{conjecture}

We next discuss the behavior of the twisted NZ matrices under a 2--3 Pachner move
relating an ideal triangulation $\calT$ with $N$ tetrahedra to another one
$\overline{\calT}$ with $N+1$ tetrahedra.  
Recall that such a move is determined by two tetrahedra with a common face in $\calT$,
which become three tetrahedra with a common edge in $\overline{\calT}$ as shown in
Figure~\ref{fig.pachner} below. Let us write the twisted NZ matrices $\mathbf{A}(t)$ and
$\mathbf{B}(t)$ of the triangulation $\calT$ schematically in columns as
\be 
\mathbf{A}(t) = (a_1,\;a_2,\;a_i)\,,\quad \mathbf{B}(t) = (b_1,\;b_2,\;b_i)\,, 
\ee
with $a_i$ meaning $(a_3,a_4,...,a_N)$ and similarly for $b_i$.
Let $\overline{\mathbf{A}}(t)$ and $\overline{\mathbf{B}}(t)$ denote
the corresponding twisted NZ matrices of $\overline{\calT}$.


\begin{theorem}
\label{thm.NZpachner}
With the above notation, there exists  a lower triangular matrix $P$ with 1's on the diagonal such that 
\begin{equation} 
\label{3AB}
P \, \overline{\mathbf{A}}(t) = \begin{pmatrix} -1 & -1 & -1  & 0 \\ 
  b_1+b_2 & a_1 & a_2 & a_i \end{pmatrix}\,,\qquad
P \, \overline{\mathbf{B}}(t) = \begin{pmatrix} -1 & -1 & -1 & 0 \\
  0 & a_2+b_1 & a_1+b_2 & b_i 
\end{pmatrix} .
\end{equation}
\end{theorem}
This corrects the omission of $P$ in~\cite[Eqn.(3.27)]{DG1}, which does not affect
the validity of the proofs in~\cite{DG1}.

\subsection{The 1-loop invariant and its twisted version}
\label{sub.1loopt}

The adjoint Reidemeister torsion has a conjectural description in terms of the
1-loop invariant of~\cite{DG1}. The latter depends on the NZ matrices
of an ideal triangulation, its shapes, their flattenings and a peripheral
curve $\gamma$. With the notation of Section~\ref{sub.tNZ}
the \emph{1-loop invariant} is defined by
\be
\label{1loop}
\begin{aligned}
  \tau^\CS_\gamma(\calT)
& := \frac{\det \left(  \widehat{\mathbf{A}}  \, \mathrm{diag}(\zeta) +
    \widehat{\mathbf{B}} \, \mathrm{diag}(\zeta'') \right)}{
2\prod_{j=1}^N\zeta_j^{f_j} \zeta_j'^{f_j'} \zeta_j''^{f_j''}} \\
& \,\,=
\frac{\det \left(  \widehat{\mathbf{G}}\, \mathrm{diag}(\zeta)
	+\widehat{\mathbf{G}}' \, \mathrm{diag}(\zeta')+\widehat{\mathbf{G}}'' \,
	\mathrm{diag}(\zeta'') \right)}{
	2\prod_{j=1}^N\zeta_j^{f_j} \zeta_j'^{f_j'} \zeta_j''^{f_j''}} 
 \in F/ (\pm1)
\end{aligned}
\ee
where the last equality follows from Equation \eqref{3GAB}.
In~\cite{DG1} it is conjectured that the 1-loop invariant $\tau^\CS_\gamma(\calT)$ is
equal to the adjoint torsion $\tau_\gamma(M)$ with respect to
$\gamma$
\be
\label{conj.DG}
\tau^\CS_\gamma(\calT) \stackrel{?}{=}  \tau_\gamma(M)  \in F^\times / (\pm1)  .
\ee

Given the above conjecture and the discussion of Section~\ref{sub.tt}, it
is natural to predict the existence of a twisted version of the 1-loop invariant,
defined in terms of the twisted NZ matrices. With the notation of
Section~\ref{sub.tNZ}, we define the \emph{twisted 1-loop invariant}  by 
\be
\label{1loopt}
\begin{aligned}
\tau^\CS(\calT,t)
&:=
\frac{\det\left(
\mathbf{A}(t)  \, \mathrm{diag}(\zeta) + \mathbf{B}(t) \, \mathrm{diag}(\zeta'')
\right)}{
\prod_{j=1}^N\zeta_j^{f_j} \zeta_j'^{f_j'} \zeta_j''^{f_j''}} \\
& := \frac{\det\left(
\mathbf{G}(t)\, \mathrm{diag}(\zeta)+\mathbf{G}'(t) \,
\mathrm{diag}(\zeta')+\mathbf{G}''(t) \,
\mathrm{diag}(\zeta'') \right)}{
\prod_{j=1}^N\zeta_j^{f_j} \zeta_j'^{f_j'} \zeta_j''^{f_j''}}
\in F[t^{\pm1}]/ (\pm t^\BZ)
\end{aligned}
\ee
where the second equality follows from the fact (analogous to Equation~\eqref{3GAB})
\be
\label{3GABt}
\mathbf{G}(t)\, \mathrm{diag}(\zeta)
+\mathbf{G}'(t) \, \mathrm{diag}(\zeta')
+ \mathbf{G}''(t) \, \mathrm{diag}(\zeta'')
= \mathbf{A}(t)  \, \mathrm{diag}(\zeta) + \mathbf{B}(t) \, \mathrm{diag}(\zeta'').
\ee

An elementary observation is that $\tau^\CS(\calT,t)$ is well-defined up to
multiplication by an element in $\pm t^\BZ$. Indeed, the computation
in~\cite[Sec.3.5]{DG1} shows that up to sign, the twisted 1-loop
invariant is independent of the choice of a combinatorial flattening. 
It also manifestly independent of the choice of a quad type of $\Delta_j$, as the
definition~\eqref{1loopt} is symmetric with respect to
$\zeta,\zeta',$ and $\zeta''$. Finally, a different choice of lifts of the
edges of $\calT$ results to left multiplication of all of the matrices
$\mathbf{G}^\square(t)$ by the same diagonal matrix $\diag(t^{c_1},t^{c_2},\dots,t^{c_N})$
for integers $c_1,c_2,\dots,c_N$, hence the observation follows. 

By its very definition of $\tau^\CS(\calT,t)$, it leads to a twisted version of
Conjecture~\eqref{conj.DG}, namely 
\be
\label{TT}
\tau^\CS(\calT,t) \stackrel{?}=  \tau(M,t) \in F[t^{\pm1}]/ (\pm t^\BZ)  .
\ee

We now list some properties of the twisted 1-loop invariant. The first concerns
the topological invariance.

\begin{theorem}
\label{thm.top}
$\tau^\CS(\calT,t)$ is invariant under Pachner 2--3 moves between ideal triangulations
that support the geometric representation.
\end{theorem}
Combining the above theorem with Proposition 1.7 of~\cite{DG1}, we conclude that the
twisted 1-loop invariant defines a topological invariant of one-cusped hyperbolic
3-manifolds.

The next two theorems concern the properties~\eqref{key1} and~\eqref{key2} of the
adjoint twisted Alexander polynomial.

\begin{theorem}
\label{thm.1}
For all $n \geq 1$ we have
\be
\label{1loopcyclic}
  \tau^{\CS}(\calT^{(n)},t^n)=\prod_{\omega^n=1}\tau^{\CS}(\calT,\omega\, t)  .
\ee
\end{theorem}

\begin{theorem}
\label{thm.2}
We have
$\tau^\CS(\calT,1)=0$ and $ \left. \frac{d}{dt}\right|_{t=1}
\tau^{\CS}(\calT,t) = \tau_\lambda^{\CS}(\calT)$
where $\lambda$ is the peripheral curve satisfying $\alpha(\lambda)=0$.
\end{theorem}

A corollary of the above theorems is a relation between the twisted 1-loop invariant
and the 1-loop invariant of~\cite{DG1} for cyclic covers, namely
\be
\label{behavior2}
\frac{\tau^\CS_\mu(\calT^{(n)})}{\tau^\CS_\mu(\calT)}
=   \prod_{ \substack{\omega^n=1 \\ \omega \neq 1}} \tau^\CS(\calT,w)  .
\ee

The next result concerns the symmetries of the twisted 1-loop invariant.
Conjecture~\ref{conj.ABt} implies a symmetry of the twisted 1-loop invariant which is
known to hold for the adjoint twisted Alexander polynomial~\cite{Kitano96,KL99},
namely $\tau(M,t)=\tau(M,1/t)$. 

\begin{corollary}(assuming Conjecture~\ref{conj.ABt})
\label{cor.sym}
We have
\be
\label{1loopsymm}
\tau^\CS(\calT,t) = \tau^\CS(\calT,1/t) \in F[t^{\pm1}]/(\pm t^\BZ)  . 
\ee
\end{corollary}

We end our discussion on the twisted 1-loop invariant with a remark which suggests
that the $t$-deformation variable is independent from the variable $\mathfrak{m}$
which is an eigenvalue of the meridian of an $\SL_2(\BC)$-representation of $\pi_1(M)$. 

\begin{remark}
\label{rem.tcharacter}  
By varying the representation, the adjoint Reidemeister torsion can be extended to
a rational function on the geometric component $X_M$ of the $\SL_2(\BC)$-character
variety of a knot complement $M$~\cite{Dubois-Garoufalidis}, and consequently
to a rational function of the geometric component of the $A$-polynomial curve.
The same holds for the 1-loop invariant; see~\cite[Sec.4]{DG1}. Likewise, we can
extend the twisted 1-loop invariant to a twisted rational function on the geometric
component of the character variety, i.e., to an element of the ring
$C(X_M)[t^{\pm1}]/(\pm t^\BZ)$, where $C(X_M)$ is the field of rational functions on
$X_M$. Then, $\mathfrak{m} \in C(X_M)$ is a variable independent of $t$.
\end{remark}



\section{Proofs: twisted NZ matrices}
\label{sec.part1}

In this section we give proofs of the properties of the twisted NZ
matrices, namely the relation to the NZ matrices of cyclic covers
(Theorem~\ref{thm.NZcyclic}), the symplectic properties (Theorem~\ref{thm.sympl}),
and the behavior under 2--3 Pachner moves (Theorem~\ref{thm.NZpachner}).

\begin{proof}(of Theorem~\ref{thm.NZcyclic})
It is clear that $\calT^{(n)}$ has $nN$ tetrahedra, as $\calT$ has $N$
tetrahedra  $\Delta_j$, $j=1,\ldots,N$. Fixing  a lift $\ti \Delta_j$ of
each $\Delta_j$ to the infinite cyclic cover, we choose lifts of
the tetrahedra of $\calT^{(n)}$ by  $\mu^k \cdot \ti \Delta_j$ for $k=0, \ldots, n-1$.
Similarly,  we choose lifts of the edges of $\calT^{(n)}$. Then 
 it is clear from the
construction of the $n$-fold cyclic cover that
\be
\label{Gnkt}
\mathbf{G}^{(n)\square}_k = 
\begin{pmatrix*}[l]
  \mathbf{G}^\square_{nk} &  \mathbf{G}^\square_{nk+1}& \cdots
  &  \mathbf{G}^\square_{nk+n-1}\\
  \mathbf{G}^\square_{nk-1} &  \mathbf{G}^\square_{nk}  & \cdots
  &  \mathbf{G}^\square_{nk+n-2}\\
\quad \vdots & \quad \vdots & \ddots & \quad \vdots \\
\mathbf{G}^\square_{nk-n+1} &  \mathbf{G}^\square_{nk-n+2}  & \cdots
&  \mathbf{G}^\square_{nk}
\end{pmatrix*} 
\ee	
for all $k \in \BZ$, which implies 
\be
\label{Gnt}
\mathbf{G}^{(n)\square} = 
\begin{pmatrix*}[c]
  \sum_{r \equiv 0} \mathbf{G}^\square_r & \sum_{r \equiv 1} \mathbf{G}^\square_r
  & \cdots & \sum_{r \equiv n-1} \mathbf{G}^\square_r\\
  \sum_{r \equiv n-1} \mathbf{G}^\square_r & \sum_{r \equiv 0} \mathbf{G}^\square_r
  & \cdots & \sum_{r \equiv n-2} \mathbf{G}^\square_r\\
  \vdots & \vdots & \ddots & \vdots \\
  \sum_{r \equiv 1} \mathbf{G}^\square_r & \sum_{r \equiv 2} \mathbf{G}^\square_r
  & \cdots & \sum_{r \equiv 0} \mathbf{G}^\square_r
\end{pmatrix*} 
\ee
where $\equiv$ means the equality of integers in modulo $n$. This 
completes the proof of Theorem~\ref{thm.NZcyclic}.
\end{proof}

\begin{proof}(of Theorem~\ref{thm.sympl})
Following Neumann~\cite{Neumann04}, we consider the $\BZ[t^{\pm1}]$-module $V$
generated by variables $Z_j,Z'_j,Z''_j$ subject to linear relations $Z_j+Z'_j+Z''_j=0$
for $ 1 \leq j \leq N$. We define a symplectic form on $V$ as
\be
\label{eqn.inner}
\langle Z_i, Z''_j \rangle = \langle Z'_i, Z_j \rangle = \langle Z''_i, Z'_j \rangle
= \delta_{ij}
\ee
and 
\[
\langle c Z_ j, Z''_j \rangle= \langle Z_ j,  \bar{c} Z''_j \rangle = c, \quad
c \in \BZ[t^{\pm1}]
\]
where $\delta_{ij}$ is the Kronecker delta and
$\bar{\ } : \BZ[t^{\pm1}]\rightarrow\BZ[t^{\pm1}]$ is the involution induced
from the inversion $t \mapsto t^{-1}$. Note that
$\langle Z''_i, Z_j \rangle = \langle Z''_i, -Z'_j - Z''_j\rangle = -\delta_{ij}$
and similarly, $\langle Z_i,Z'_j\rangle = \langle Z'_i,Z''_j \rangle = - \delta_{ij}$.
	
We associate an element $E_i \in V$ to each edge $e_i$ of $\calT$ as
\be
\label{eqn.ei}
E_i =  \sum_{j=1}^N  \mathbf{G}(t)_{ij}Z_j+ \mathbf{G}'(t)_{ij}Z'_j
+\mathbf{G}''(t)_{ij}Z''_j  .
\ee
In what follows, we use the same argument given in the proof
of~\cite[Thm.3.6]{choi06} to show that
\be
\label{inner}
\langle E_a, E_b \rangle=0  \quad \textrm{for all } 1 \leq a,b \leq N  .
\ee
Considering the linear expansion of $\langle E_a,E_b \rangle$ using the
equation~\eqref{eqn.ei}, we have a non-trivial term  when the edges $e_a$ and $e_b$
appear in some $\Delta_j$ with different shape parameters. Without loss of generality,
suppose that $e_a$ and $e_b$ have parameters $Z_j$ and $Z''_j$, respectively.
It is clear that  $\widetilde{\Delta}_j$ contains $\mu^\alpha \cdot \widetilde{e}_a$ and 
$\mu^\beta \cdot \widetilde{e}_b$ for some $\alpha$ and $\beta \in \BZ$. Equivalently,
$\mu^{-\alpha} \cdot \widetilde{\Delta}_j$ and $\mu^{-\beta} \cdot \widetilde{\Delta}_j$
are attached to $\widetilde{e}_a$  and $\widetilde{e}_b$, respectively. We thus get a
non-trivial term 
\be 
\label{eqn.contr}
\langle t^{-\alpha} Z_j, t^{-\beta} Z''_j \rangle= t^{\beta-\alpha} 
\ee
in the linear expansion from the triple $(\widetilde{\Delta}_j,
\mu^\alpha \cdot \widetilde{e}_a, \mu^\beta \cdot \widetilde{e}_b)$. On the other hand,
considering the face-pairing of the face of $\Delta_j$ containing $e_a$ and $e_b$,
there is another tetrahedron  $\Delta_{i}$ (possibly $i=j$), which also contains both
$e_a$ and $e_b$. It follows that $\widetilde{\Delta}_i$ contains
$\mu^{\alpha+\gamma} \cdot \widetilde{e}_a$ and $\mu^{\beta+\gamma} \cdot \widetilde{e}_b$
for some $\gamma \in \BZ$. From this triple $(\widetilde{\Delta}_i,
\mu^{\alpha+\gamma} \cdot \widetilde{e}_a, \mu^{\beta+\gamma} \cdot \widetilde{e}_b)$, we
get another non-trivial term $ -t^{(\beta+\gamma)-(\alpha+\gamma)}=-t^{\beta-\alpha}$ which
we can pair with the term~\eqref{eqn.contr} to cancel out. Note that we have the minus
sign here, since the face-pairing is orientation-reversing: the shape parameters of $e_a$
and $e_b$ for $\Delta_i$ are either $(Z_i,Z'_i)$, $(Z'_i,Z''_i)$, or $(Z''_i,Z_i)$.
In this way, every non-trivial term that appears in the expansion of
$\langle E_a,E_b\rangle$ is paired to cancel out. This proves the equation~\eqref{inner}. 
We then directly obtain the desired  equation~\eqref{eqn.abba}, since we have 
\begin{align*}
& \langle E_a, E_b \rangle  =
\sum_{j=1}^N \left\langle \mathbf{G}(t)_{aj} Z_j+\mathbf{G}'(t)_{aj}
Z'_j+\mathbf{G}''(t)_{aj} Z''_j,
\mathbf{G}(t)_{bj} Z_j+\mathbf{G}'(t)_{bj} Z'_j+\mathbf{G}''(t)_{bj} Z''_j\right
\rangle
\\
&=\sum_{j=1}^N \mathbf{G}(t)_{aj} (\mathbf{G}''(1/t)_{bj}
-\mathbf{G}'(1/t)_{bj})+\mathbf{G}'(t)_{aj} (\mathbf{G}(1/t)_{bj}
-\mathbf{G''}(1/t)_{bj})
+\mathbf{G}''(t)_{aj} (\mathbf{G}'(1/t)_{bj}-\mathbf{G}(1/t)_{bj})
\\
&=\sum_{j=1}^N (\mathbf{G}(t)_{aj}-\mathbf{G}'(t)_{aj})
(\mathbf{G}''(1/t)_{bj}-\mathbf{G}'(1/t)_{bj})-
(\mathbf{G}''(t)_{aj}-\mathbf{G}'(t)_{aj})
(\mathbf{G}(1/t)_{bj}-\mathbf{G}'(1/t)_{bj})
\\
&= \sum_{j=1}^N \mathbf{A}(t)_{aj} \mathbf{B}(1/t)_{bj}
- \mathbf{B}(t)_{aj} \mathbf{A}(1/t)_{bj}  .
\end{align*}
\end{proof}

\begin{proof}(of Theorem~\ref{thm.NZpachner})
Recall that a 2--3 Pachner move on $\calT$ is determined by two tetrahedra
$\Delta_\alpha$ and $\Delta_\beta$ with a common face in $\calT$. Dividing the
bipyramid $\Delta_\alpha \cup \Delta_\beta$ into three tetrahedra $\Delta_a, \Delta_b,$
and $\Delta_c$ as in Figure~\ref{fig.pachner}, we obtain a new ideal triangulation
$\overline{\calT}$ with one additional edge $e_0$.
We choose lifts of the five tetrahedra $\Delta_\alpha$, $\Delta_\beta$, $\Delta_a$,
$\Delta_b$, $\Delta_c$ and the edge $e_0$ (to the infinite cyclic cover)
such that $\ti \Delta_\alpha \cup \ti \Delta_\beta$ and
$\ti \Delta_a \cup \ti \Delta_b \cup \ti \Delta_c$ are the same bipyramid
containing $\ti e_0$. Note that 
the (twisted) gluing equation of $e_0$ is 
\be
\label{eqn.g0}
\log z'_a + \log z'_b + \log z'_c = 2 \pi \sqrt{-1}.
\ee

\begin{figure}[htpb!]
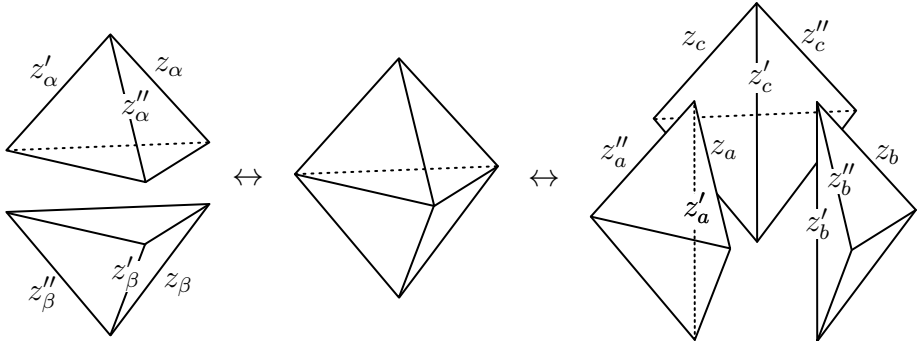
\caption{2-3 Pachner move.}
\label{fig.pachner}
\end{figure}

There are nine edges on the bipyramid whose gluing equations are affected by the
Pachner move: three at the triangle of the base pyramid and six others
on the bipyramid but not on its base. The six edges give
\be
\label{eqn.relation}
\begin{aligned}
\log z_\alpha & =\log z_b + \log z_c'' &
\log z_\alpha' &=\log z_a''+\log z_c   &
\log z_\alpha'' & =\log z_a+\log z_b''
\\
\log z_\beta & =\log z_b''+\log z_c &
\log z_\beta' & = \log z_a'' +\log z_b &
\log z_\beta''&= \log z_a+ \log z_c'',
\end{aligned}
\ee	
and the three edges give
\be
\label{eqn.relation2}
\begin{aligned}
  \log z_\alpha + \log z_\beta & =\log z'_a &
  \log z'_\alpha + \log z''_\beta & =\log z'_b &
  \log z''_\alpha + \log z'_\beta=\log z'_c.
\end{aligned}
\ee	
The relations~\eqref{eqn.relation2} are obtained from the relations~\eqref{eqn.relation}
up to the gluing equation \eqref{eqn.g0} of $e_0$ and the relation
$\log z + \log z' + \log z''= \pi \sqrt{-1}$. For instance,
\begin{align*}
	\log z_\alpha + \log z_\beta &= \log z_b + \log z''_c + \log z''_b +\log z_c \\
	& = 2 \pi \sqrt{-1} - \log z'_b - \log z'_c \\
	&= \log z'_a  .
\end{align*}
Therefore, we can obtain the twisted gluing equation matrices
$\overline{\mathbf{G}}^\square(t)$ for $\overline{\calT}$ from the
twisted gluing equation matrices for $\overline{\calT}$ by substituting
$\log z_\alpha^\square$ and $\log z_\beta^\square$ in terms of the right hand side
of equation~\eqref{eqn.relation} and the substitution is correct up to the
gluing equation of $e_0$ and the relation $\log z + \log z' + \log z''= \pi \sqrt{-1}$.
In matrix form, this implies that for some lower triangular matrix
$Q=\begin{pmatrix} 1 & 0 \\ * & I_N \end{pmatrix}$ and an integer matrix $C$, we have
{\small
\begin{equation}
\label{pre5G}  
\begin{aligned}
\begin{pmatrix}
\overline{\mathbf{G}}(t)_{a} &
\overline{\mathbf{G}}(t)_{b} & 
\overline{\mathbf{G}}(t)_{c}
\end{pmatrix}
&= Q\begin{pmatrix}
0& 0 & 0\\
\mathbf{G}''(t)_{\alpha} + \mathbf{G}''(t)_{\beta} &  
\mathbf{G}(t)_{ \alpha} + \mathbf{G}'(t)_{\beta}& 
\mathbf{G}'(t)_{\alpha} + \mathbf{G}(t)_{\beta} 
\end{pmatrix} + C, \\
\begin{pmatrix}
\overline{\mathbf{G}}'(t)_{a} &
\overline{\mathbf{G}}'(t)_{b} & 
\overline{\mathbf{G}}'(t)_{c} 
\end{pmatrix}
&= Q\begin{pmatrix}
1& 1 & 1\\
\mathbf{0} &  \mathbf{0}& \mathbf{0}
\end{pmatrix} + C, \\
\begin{pmatrix}
\overline{\mathbf{G}}''(t)_{a} &
\overline{\mathbf{G}}''(t)_{b} & 
\overline{\mathbf{G}}''(t)_{c} 
\end{pmatrix}
&= Q \begin{pmatrix}
0& 0 & 0\\
\mathbf{G}'(t)_{\alpha} + \mathbf{G}'(t)_{ \beta} &  
\mathbf{G}''(t)_{ \alpha} + \mathbf{G}(t)_{ \beta}& 
\mathbf{G}(t)_{\alpha} + \mathbf{G}''(t)_{\beta} 
\end{pmatrix} + C  .
\end{aligned}
\end{equation}
}
Here $I_N$ is the identity matrix of size $N$ and $X_{j}$ is the $j$-th column of
a matrix $X$. Multiplying  $P=\begin{pmatrix} 1 & 0 \\
\mathbf{G}'(t)_{\alpha} + \mathbf{G}'(t)_{\beta} & I_N \end{pmatrix} Q^{-1}$ on both
sides of Equation \eqref{pre5G}, we obtain
{\small
\begin{equation}
\label{5G}  
\begin{aligned}
P\begin{pmatrix}
\overline{\mathbf{G}}(t)_{a} &
\overline{\mathbf{G}}(t)_{b} & 
\overline{\mathbf{G}}(t)_{c}
\end{pmatrix}
&= \begin{pmatrix}
0& 0 & 0\\
\mathbf{G}''(t)_{\alpha} + \mathbf{G}''(t)_{\beta} &  
\mathbf{G}(t)_{ \alpha} + \mathbf{G}'(t)_{\beta}& 
\mathbf{G}'(t)_{\alpha} + \mathbf{G}(t)_{\beta} 
\end{pmatrix} + PC, \\
P  \begin{pmatrix}
\overline{\mathbf{G}}'(t)_{a} &
\overline{\mathbf{G}}'(t)_{b} & 
\overline{\mathbf{G}}'(t)_{c} 
\end{pmatrix}&=\begin{pmatrix}
1& 1 & 1\\
\mathbf{G}'(t)_{ \alpha} + \mathbf{G}'(t)_{\beta} &  
\mathbf{G}'(t)_{ \alpha} + \mathbf{G}'(t)_{\beta} & 
\mathbf{G}'(t)_{ \alpha} + \mathbf{G}'(t)_{\beta}
\end{pmatrix}+ PC, \\
P \begin{pmatrix}
\overline{\mathbf{G}}''(t)_{a} &
\overline{\mathbf{G}}''(t)_{b} & 
\overline{\mathbf{G}}''(t)_{c} 
\end{pmatrix}
&= \begin{pmatrix}
0& 0 & 0\\
\mathbf{G}'(t)_{\alpha} + \mathbf{G}'(t)_{\beta} &  
\mathbf{G}''(t)_{\alpha} + \mathbf{G}(t)_{\beta}& 
\mathbf{G}(t)_{\alpha} + \mathbf{G}''(t)_{\beta} 
\end{pmatrix} + PC \, .
\end{aligned}
\end{equation}
}
It follows that  
\begin{equation*}
\begin{aligned}
P\begin{pmatrix}
\overline{\mathbf{A}}(t)_{a} &
\overline{\mathbf{A}}(t)_{b} & 
\overline{\mathbf{A}}(t)_{c}
\end{pmatrix}
&= \begin{pmatrix}
-1& -1 & -1\\
\mathbf{B}(t)_{ \alpha} + \mathbf{B}(t)_{\beta} &  
\mathbf{A}(t)_{ \alpha} & 
\mathbf{A}(t)_{\beta} 
\end{pmatrix}, \\
P \begin{pmatrix}
\overline{\mathbf{B}}(t)_{a} &
\overline{\mathbf{B}}(t)_{b} & 
\overline{\mathbf{B}}(t)_{c} 
\end{pmatrix}
&= \begin{pmatrix}
-1& -1 & -1\\
\mathbf{0} &  
\mathbf{A}(t)_{ \beta}+ \mathbf{B}(t)_{ \alpha} & 
\mathbf{A}(t)_{ \alpha} + \mathbf{B}(t)_{\beta} 
\end{pmatrix},
\end{aligned}
\end{equation*}
where $P$ is a lower triangular matrix with 1's on the diagonal.
\end{proof}


\section{Examples}
\label{sec.examples}

In this section we explain how to compute the twisted NZ matrices
of an ideal triangulation using the methods of \texttt{SnapPy}.
As is customary, we do so by the example of the simplest hyperbolic $4_1$
knot, and later by the $6_3$ knot. 

\subsection{The  knot $4_1$}
\label{sub.41a}

\begin{figure}[htpb!]
	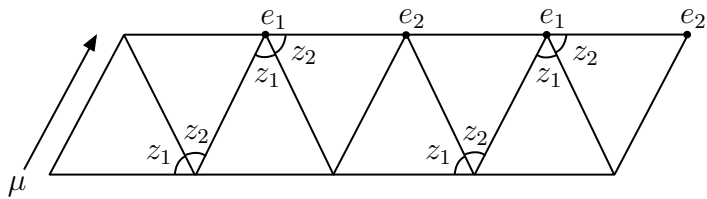
	\caption{The triangulation of the peripheral torus of $4_1$.}
	\label{fig.torus41}
\end{figure}

The knot complement $M$ of the knot $4_1$ decomposes into two ideal tetrahedra
$\Delta_1$ and $\Delta_2$. Using the default triangulation of \texttt{SnapPy}
for $4_1$ with isometry signature \texttt{cPcbbbiht\_BaCB}, 
the gluing equations are
\begin{align*}
	e_1&: z_1^2 z_1' z_2^2 z_2'=1, \\
	e_2&: z_1' (z_1'')^2 z_2' (z_2'')^2=1 \,, 
\end{align*}
the corresponding gluing equation matrices are
\[
\mathbf{G}
= \begin{pmatrix}
  2 & 2 \\
  0 & 0
\end{pmatrix},
\qquad
\mathbf{G}'
= \begin{pmatrix}
  1 & 1 \\
  1 & 1
\end{pmatrix}, \qquad
\mathbf{G}''
= \begin{pmatrix}
  0 & 0 \\
  2 & 2
\end{pmatrix} \,,
\]
and the NZ matrices are
\[
\mathbf{A}= \begin{pmatrix} 1  &  1\\ -1 & -1 \end{pmatrix}, \qquad
\mathbf{B}= \begin{pmatrix} -1  & - 1 \\ 1 & 1 \end{pmatrix} 
\]
with
\[
\mathbf{A} \mathbf{B}^T = \begin{pmatrix} -2 & 2 \\ 2 & -2 \end{pmatrix}   
\]
symmetric. We can read the above information by looking at the triangulated cusp
shown in Figure~\ref{fig.torus41}.

To describe the ideal triangulation $\ti \calT$ of the infinite cyclic 
cover $\ti M$ of $M$, we consider the triangulation of the peripheral torus induced
from $\calT=\{\Delta_1,\Delta_2\}$ (see Figure \ref{fig.torus41}) together with its
lift to the universal cover, the Euclidean plane $\BE^2$. Note that each triangle
of $\BE^2$ with corner $z_j$ has a label $j \in \{1,2\}$ of the corresponding
tetrahedron. 

\begin{figure}[htpb!]
	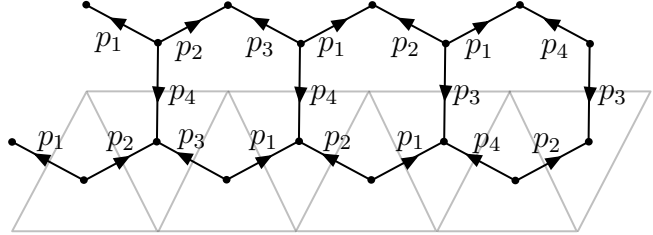
	\caption{Face-pairings}
	\label{fig.example41}
\end{figure}

We now describe a second labeling of each triangle of $\BE^2$ by an integer,
obtained by choosing a cocycle representative of the abelianization map
$\alpha: \pi_1(M) \twoheadrightarrow \BZ$. This construction is standard and	
is appears in unpublished work of Goerner and also in
Zickert~\cite[Sec.3.5]{Zickert:Apoly}. Let $G$ denote the dual dual 1-skeleton
of the triangulation of $\BE^2$ (see Figure~\ref{fig.example41}). Each vertex and
edge of $G$ correspond to a tetrahedron and face-pairing of $\calT$, respectively,
and that the
edges of $G$ are oriented from $\Delta_0$ to $\Delta_1$. Since the face-pairings
generate the fundamental group $\pi_1(M)$, one can assign an integer $\alpha(p_i)$ to
each face-pairing $p_i$ corresponding to the abelianization map
$\alpha : \pi_1(M) \twoheadrightarrow \BZ$. Since in our case there are two tetrahedra,
we have four face-pairings, thus four integers $\alpha(p_1),\ldots,\alpha(p_4)$
that are required to satisfy
\begin{equation}
\label{ap41}
\begin{aligned}
	\alpha(p_2)-\alpha(p_4) & =1  &	\alpha(p_1)-\alpha(p_4)& =1,\\
	\alpha(p_1)-\alpha(p_3) & =1  & \alpha(p_2)-\alpha(p_3)& =1  .
\end{aligned}
\end{equation}      
Note that these equations are directly obtained from Figure \ref{fig.example41} with
the $\mu$-action described in Figure~\ref{fig.torus41}.

\begin{figure}[htpb!]
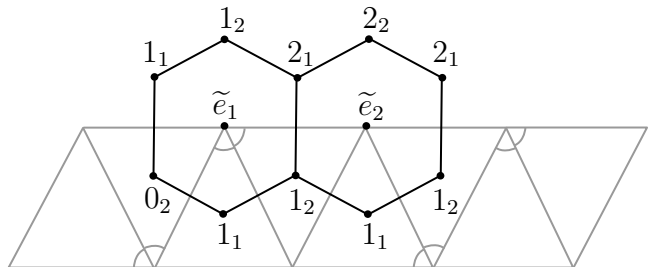
\caption{The dual 1-complex of $\BE^2$, where the labeling $\alpha(v)_j$ of a vertex
  $v$ indicates that $v$ is labeled by the $j$-th tetrahedron.}
\label{fig.example41v}
\end{figure}

Since $\BE^2$ is a contractible space, we can assign an integer $\alpha(v)$ to each
vertex $v$ of $G$ such that $\alpha(v_2)-\alpha(v_1)=\alpha([v_1,v_2])$, where
$[v_1,v_2]$ is an oriented edge of $G$ from $v_1$ to $v_2$. Explicitly,
the general solution of the linear equations~\eqref{ap41} is given by
\be
(\alpha(p_1), \alpha(p_2), \alpha(p_3), \alpha(p_4))=(a,a,a-1,a-1) \qquad
(a \in \BZ)  .
\ee
It follows from the construction that if a vertex $v$ corresponds to $\Delta_j$
in $\calT$, then it corresponds to $\mu^{\alpha(v)} \cdot \widetilde{\Delta}_j$
in $\ti \calT$. Choosing $(\alpha(p_1),\alpha(p_2),\alpha(p_3),\alpha(p_4))=(0,0,-1,-1)$
and $\alpha(v)$ and $\widetilde{e}_i$ as in Figure~\ref{fig.example41}, we obtain
the table
\begin{center}	
\begin{tabular}{c|c|c|c}
& 0 & 1 & 2\\
\hline
$\widetilde{e}_1$ & $z'_2$ & $z_1^2z_2^2 $& $z_1'$ \\
$\widetilde{e}_2$ &  & $z_1' (z_2'')^2$ & $(z_1'')^2 z_2'$
\end{tabular}
\end{center}
where the $(i,k)$-entry is the product of shape-parameters contributed from
$\mu^k \cdot \widetilde{\Delta}_1$ and $\mu^k \cdot \widetilde{\Delta}_2$ to
$\widetilde{e}_i$ (the empty entries are $1$). The gluing equation around
$\widetilde{e}_i$ is the product of the entries of the corresponding row in the
above table, i.e., the product over the vertices of polygons around each edge
$\widetilde{e}_i$ in Figure~\ref{fig.example41}. Hence, we have 
\be
\label{Gt41}
\mathbf{G}(t) = \begin{pmatrix} 2 t  &  2t \\ 0 & 0 \end{pmatrix}, \qquad
\mathbf{G}'(t) = \begin{pmatrix} t^2  &  1 \\ t & t^2 \end{pmatrix}, \qquad
\mathbf{G}''(t) = \begin{pmatrix} 0  &  0 \\ 2t^2 & 2t \end{pmatrix},
\ee
\be
\label{ABt41}
\mathbf{A}(t) = \begin{pmatrix} - t^2 + 2 t &  2t -1\\ -t & -t^2 \end{pmatrix}, \qquad
\mathbf{B}(t) = \begin{pmatrix} -t^2  & - 1 \\ 2t^2-t & -t^2+2t \end{pmatrix} \,, 
\ee
giving that
\[
\mathbf{A}(t)\mathbf{B}(1/t)^T =
\begin{pmatrix}
-2t+2-2 t^{-1} & t + t^{-2} \\
t^{2} + t^{-1} & -2t+2-2 t^{-1}
\end{pmatrix} \,,
\]
which specializes to the matrix $\mathbf{A} \mathbf{B}^T$ when $t=1$.

\subsection{The  knot $6_3$}
\label{sub.63a}

\begin{figure}[htpb!]
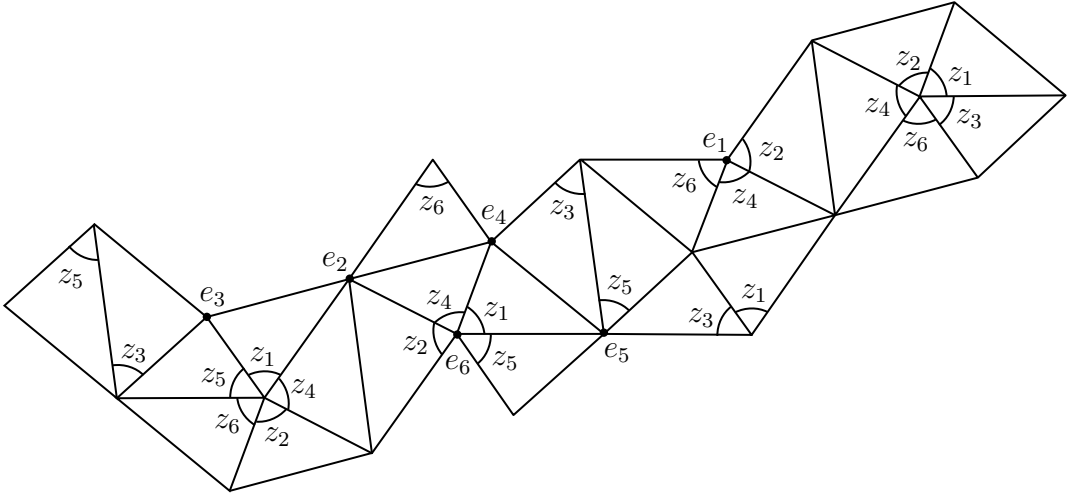
\caption{The triangulation of the peripheral torus of $6_3$.}
\label{fig.torus63}
\end{figure}

The default \texttt{SnapPy} triangulation of the knot $6_3$ with isometry signature
\texttt{gLLPQccdefffhggaacv\_aBBb} has six tetrahedra $\Delta_1,\ldots,\Delta_6$
with the peripheral torus shown in Figure~\ref{fig.torus63}. 

The gluing equations of the six edges are 
\begin{align*}
e_1 & : z_1 z_2 z_3 z_4 z_6 =1,\\
e_2 & : z'_1 z_2' (z''_2)^2 z'_4 (z''_4)^2 z'_6=1,\\
e_3 & : z''_1 z'_2 z'_3 z'_5 z''_6 =1, \\
e_4 & : z''_1 z'_3 z'_4 z'_5 z''_6 =1, \\
e_5 & : z'_1 z_3 (z''_3)^2 z_5 (z''_5)^2 z'_6 = 1, \\
e_6 &:  z_1 z_2 z_4 z_5 z_6 =1 
\end{align*}
and the corresponding gluing equation matrices are 
{\tiny
\be
\mathbf{G} = \begin{pmatrix}
1 & 1 & 1 & 1 & 0 & 1 \\
0 & 0 & 0 & 0 & 0 & 0 \\
0 & 0 & 0 & 0 & 0 & 0 \\
0 & 0 & 0 & 0 & 0 & 0 \\
0 & 0 & 1 & 0 & 1 & 0 	\\
1 & 1 & 0 & 1 & 1 & 1
\end{pmatrix}, \qquad
\mathbf{G}' = \begin{pmatrix}
0 & 0 & 0 & 0 & 0 & 0 \\
1 & 1 & 0 & 1 & 0 & 1 \\
0 & 1 & 1 & 0 & 1 & 0 \\
0 & 0 & 1 & 1 & 1 & 0 \\
1 & 0 & 0 & 0 & 0 & 1 	\\
0 & 0 & 0 & 0 & 0 & 0
\end{pmatrix}, \qquad
\mathbf{G}'' = \begin{pmatrix}
0 & 0 & 0 & 0 & 0 & 0 \\
0 & 2 & 0 & 2 & 0 & 0 \\
1 & 0 & 0 & 0 & 0 & 1 \\
1 & 0 & 0 & 0 & 0 & 1 \\
0 & 0 & 2 & 0 & 2 & 0 	\\
0 & 0 & 0 & 0 & 0 & 0
\end{pmatrix}  .
\ee
}

\begin{figure}[htpb!]
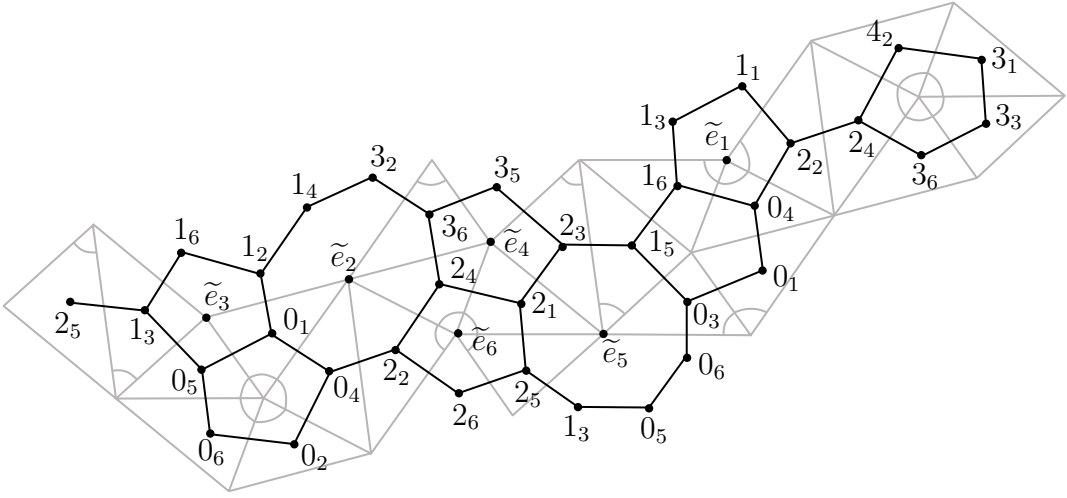
\caption{The dual 1-skeleton of $\BE^2$.}
\label{fig.example63}
\end{figure}

\noindent
The NZ matrices are
{\tiny
\[
\mathbf{A}= \begin{pmatrix}
1 & 1 & 1 & 1 & 0 & 1 \\ -1 & -1 & 0 & -1 & 0 & -1 \\ 0 & -1 & -1 & 0 & -1 & 
  0 \\ 0 & 0 & -1 & -1 & -1 & 0 \\ -1 & 0 & 1 & 0 & 1 & -1 \\ 1 & 1 & 0 & 1 & 1 & 1
\end{pmatrix}, \qquad
  \mathbf{B}= \begin{pmatrix}
0 & 0 & 0 & 0 & 0 & 0 \\ -1 & 1 & 0 & 1 & 0 & -1 \\ 1 & -1 & -1 & 0 & -1 & 1 \\ 1 & 
  0 & -1 & -1 & -1 & 1 \\ -1 & 0 & 2 & 0 & 2 & -1 \\ 0 & 0 & 0 & 0 & 0 & 0
\end{pmatrix}, 
\]
}
with
{\tiny
\[
  \mathbf{A} \mathbf{B}^T = \begin{pmatrix}
0 & 0 & 0 & 0 & 0 & 0 \\ 0 & 0 & -1 & -1 & 2 & 0 \\ 0 & -1 & 3 & 2 & -4 & 
  0 \\ 0 & -1 & 2 & 3 & -4 & 0 \\ 0 & 2 & -4 & -4 & 6 & 0 \\ 0 & 0 & 0 & 0 & 0 & 0
\end{pmatrix}   
\]
}
symmetric. As in Section~\ref{sub.41a}, we assign an integer to (each oriented edge
and) each vertex of the dual complex according to the abelianization map
$\alpha : \pi_1(M)\twoheadrightarrow \BZ$. See Figure~\ref{fig.example63}.

Choosing a lift $\widetilde{e}_i$ of $e_i$ as in Figure~\ref{fig.example63}, we have
the table
\begin{center}	
\begin{tabular}{c|c|c|c|c}
& 0 & 1 & 2 & 3 \\
\hline
$\widetilde{e}_1$ & $z_4$ & $z_1 z_3 z_6 $& $z_2$ & \\
$\widetilde{e}_2$ & $z'_1z''_4$ & $z''_2 z'_4$ & $z'_2 z''_4$ & $z''_2 z'_6$ \\
$\widetilde{e}_3$ & $z''_1 z'_5$ & $z'_2 z'_3 z''_6$ &  &  \\
$\widetilde{e}_4$ &  & & $z''_1 z'_3 z'_4$ & $z'_5 z''_6$ \\
$\widetilde{e}_5$ & $z''_3 z''_5 z'_6$ & $z_3 z_5$ & $z'_1 z''_3 z''_5$ &  \\
$\widetilde{e}_6$ &  &  & $z_1 z_2 z_4 z_5 z_6$ &  
\end{tabular}
\end{center}
where the $(i,k)$-entry is the product of shape-parameters contributed from
$\mu^k \cdot \widetilde{\Delta}_1, \ldots, \mu^k \cdot \widetilde{\Delta}_6$ to
$\widetilde{e}_i$. Hence, we have
{\tiny
\begin{align}
\label{Gt63}
\mathbf{G}(t) &= \begin{pmatrix}
t& t^2 & t & 1 & 0 & t \\
0 & 0 & 0 & 0 & 0 & 0 \\
0 & 0 & 0 & 0 & 0 & 0 \\
0 & 0 & 0 & 0 & 0 & 0 \\
0 & 0 & t & 0 & t & 0 	\\
t^2 & t^2 & 0 & t^2 & t^2 & t^2	
\end{pmatrix},
\mathbf{G}'(t) = \begin{pmatrix}
0 & 0 & 0 & 0 & 0 & 0 \\
1 & t^2 & 0 & t & 0 & t^3 \\
0 & t& t & 0 & 1 & 0 \\
0 & 0 & t^2 & t^2 & t^3 & 0 \\
t^2 & 0 & 0 & 0 & 0 & 1 	\\
0 & 0 & 0 & 0 & 0 & 0
\end{pmatrix}, 
\mathbf{G}''(t) = \begin{pmatrix}
0 & 0 & 0 & 0 & 0 & 0 \\
0 & t+t^3 & 0 & 1+t^2 & 0 & 0 \\
1 & 0 & 0 & 0 & 0 & t \\
t^2& 0 & 0 & 0 & 0 & t^3 \\
0 & 0 & 1+t^2 & 0 & 1+t^2 & 0 	\\
0 & 0 & 0 & 0 & 0 & 0
\end{pmatrix}
\end{align}
}
and
{\tiny
\begin{align}
\label{ABt63}  
\mathbf{A}(t)&= \begin{pmatrix}
t& t^2 & t & 1 & 0 & t \\
-1 & -t^2 & 0 & 0 & 0 & -t^3 \\
0 & -t & -t & 0 & -1 & 0 \\
0 & 0 & -t^2 & -t^2 & -t^3 & 0 \\
-t^2  & 0 & t & 0 & t & -1 	\\
t^2 & t^2 & 0 & t^2 & t^2 & t^2
\end{pmatrix},
\qquad
\mathbf{B}(t) = \begin{pmatrix}
0 & 0 & 0 & 0 & 0 & 0 \\
-1 & t-t^2+t^3  & 0 & 1-t+t^2 & 0 & -t^3 \\
1 & -t & -t & 0 & -1 & t \\
t^2& 0 & -t^2 & -t^2 & -t^3 & t^3 \\
-t^2 & 0 & 1+t^2 & 0 & 1+t^2 & -1 	\\
0 & 0 & 0 & 0 & 0 & 0
\end{pmatrix}, 
\end{align}
}
giving that
{\tiny
\[
\mathbf{A}(t)\mathbf{B}(1/t)^T =
\begin{pmatrix}
  0 & 0 & 0 & 0 & 0 & 0 \\
  0 & -2t+4-2 t^{-1} & -t^2+t- 1 & -1+t^{-1}- t^{-2} & t^{3}+t^{-2} & 0 \\
  0 & -1+t^{-1}- t^{-2} & 3 & t^{-1}+ t^{-3} & -t-1 -t^{-1}-t^{-2} & 0 \\
  0 & -t^2+t-1 & t^3 + t & 3 & -t^{3}-t^{2} -t^{-1}-1 & 0 \\
  0 & t^{2}+t^{-3} & -t^{2}-t^{1} -1-t^{-1}
  & -1-t^{-1} -t^{-2}-t^{-3} & 2t+2+2 t^{-1}& 0 \\
  0 & 0 & 0 & 0 & 0 & 0 
\end{pmatrix}
\]
}
which specializes to matrix $\mathbf{A} \mathbf{B}^T$ when $t=1$.


\section{Proofs: twisted 1-loop invariant}
\label{sec.part2}

In this section we give proofs of the properties of the twisted 1-loop invariant.

\begin{proof}(of Theorem~\ref{thm.top})
We will use the notation of Theorem~\ref{thm.NZpachner} which relates the NZ matrices of
two triangulations $\calT$ and $\overline{\calT}$ obtained by a 2--3 Pachner move
shown in Figure~\ref{fig.pachner}. 
Recall that a 2--3 Pachner move on $\calT$ is determined by two tetrahedra
$\Delta_\alpha$ and $\Delta_\beta$ with a common face in $\calT$. Dividing the
bipyramid $\Delta_\alpha \cup \Delta_\beta$ into three tetrahedra $\Delta_a, \Delta_b,$
and $\Delta_c$ as in Figure~\ref{fig.pachner}, we obtain a new ideal triangulation
$\overline{\calT}$ with one additional edge $e_0$.
Note that
one can read off relations of shape parameters between $\calT$ and
$\overline{\calT}$ directly from Figure~\ref{fig.pachner}:
\be
\label{eqn.shapes}
z'_a= z_\alpha z_\beta, \ z'_b=\frac{1-1/z_\beta}{1-z_\alpha},\,
z'_c = \frac{1-1/z_\alpha}{1-z_\beta} \textrm{ and }
z_\alpha = \frac{1-1/z'_b}{1-z'_c}, \,	z_\beta = \frac{1-1/z'_c}{1-z'_b}
\ee
with the gluing equation $z'_a z'_b z'_c=1$ for the edge $e_0$.

We first consider the flattening part of the definition~\eqref{1loopt}.
A combinatorial flattening of $\overline{\calT}$ satisfies the equation
$f'_a+f'_b+f'_c=2$ (coming from the edge $e_0$) and determines that  of $\calT$ by 
\be
\label{eqn.flat}
\begin{aligned}
f_\alpha & =f_b+f''_c, & f'_\alpha & = f''_a+ f_c & f''_\alpha & =f_a+f''_b,\\
f_\beta & =f''_b+f_c, &  f'_\beta & = f''_a + f_b, & f''_\beta & =f_a+f''_c.
\end{aligned}
\ee	
Then a straightforward computation shows that
\be
\label{eqn.com}
\prod_\square \zeta^{\square f^\square_a}_a \zeta_b^{\square f^\square_b}
\zeta^{\square f^\square_c}_c =
\pm \frac{z_\alpha  (1-z_\alpha)z_\beta(1-z_\beta)}{1-z_\alpha z_\beta}
\prod_\square \zeta^{\square f^\square_\alpha}_\alpha
\zeta^{\square f^\square_\beta}_\beta  .
\ee
Moreover, combining Equation \eqref{pre5G} with \eqref{eqn.shapes}, one checks that
\begin{align*}
&	
\begin{pmatrix}
\displaystyle\sum_\square \overline{\mathbf{G}}^\square(t)_{a}\zeta^\square_a  & 
\displaystyle\sum_\square \overline{\mathbf{G}}^\square(t)_{b}\zeta^\square_b & 
\displaystyle\sum_\square \overline{\mathbf{G}}^\square(t)_{c}\zeta^\square_c \end{pmatrix}   
\begin{pmatrix}
\frac{1}{z_\alpha z_\beta} & \frac{1}{z_\alpha z_\beta}-1 & \frac{1}{z_\alpha z_\beta}-1 \\
0 & \frac{z_\beta(1-z_\alpha)(1-z_\alpha z_\beta)}{1-z_\beta} & 0 \\
0 & 0 & \frac{z_\alpha(1-z_\beta)(1-z_\alpha z_\beta)}{1-z_\alpha}
\end{pmatrix}\\
& = Q \begin{pmatrix}
1& 0  & 0 \\
\ast & \displaystyle\sum_\square \mathbf{G}^\square(t)_{\alpha} \zeta^\square_\alpha
&  \displaystyle\sum_\square \mathbf{G}^\square(t)_{\beta}\zeta^\square_\beta 
\end{pmatrix}
\begin{pmatrix} 1 & 0 & 0\\
0 & z_\alpha z_\beta(1-z_\alpha) & z_\alpha(1-z_\alpha) \\
0 & z_\beta (1-z_\beta) & z_\alpha z_\beta(1-z_\beta)			
\end{pmatrix}\,,
\end{align*}
where the determinant of the $3 \times 3$ matrix in the left (resp., right) hand
side is $(1-z_\alpha z_\beta)^2$ (resp., $z_\alpha (z_\alpha-1) z_\beta (z_\beta-1)
(1-z_\alpha z_\beta))$. Combining this fact with the equation~\eqref{eqn.com} and the fact that $\det Q=1$, we
conclude that $\tau^{\CS}(\overline{\calT},t) = \pm \tau^{\CS}(\calT,t)$.
\end{proof}

\begin{proof}(of Theorem~\ref{thm.1})
Using Equation~\eqref{Gnkt} and the fact that a combinatorial flattening of
$\calT^{(n)}$ is simply given by $n$ copies of that of $\calT$, it follows that  
\begin{align}
\tau^\CS(\calT^{(n)},t^n) &= \det\left( \displaystyle\sum_{k \in \BZ}
  \left(  \sum_\square \mathbf{G}_k^{(n)\square}\,
    \mathrm{diag}(\zeta^\square,\ldots,\zeta^\square)\right) t^{nk}\right)  \cdot
\left(\prod_{j=1}^N\zeta_j^{f_j} \zeta_j'^{f_j'} \zeta_j''^{f_j''}\right)^{-n}  \\
&=\det \begin{pmatrix*}[c]
W_0 &  W_1& \cdots & W_{n-2}& W_{n-1} \\
t^{n} W_{n-1} &  W_0& \cdots & W_{n-3}& W_{n-2}  \\
\vdots & \vdots & \ddots & \vdots \\
t^{n} W_1 & t^n W_2& \cdots &t^n W_{n-1}& W_0 
\end{pmatrix*}
\cdot  \left(\prod_{j=1}^N\zeta_j^{f_j} \zeta_j'^{f_j'} \zeta_j''^{f_j''}\right)^{-n}
\label{eqn:cyctor}
\end{align}
where  $W_r = \sum_{k \in \BZ} \sum_\square \mathbf{G}_{nk+r}^{\square} \,
\mathrm{diag} (\zeta^\square)\,t^{nk} $ for  $0 \leq r \leq n-1$.
	
The block matrix in the equation~\eqref{eqn:cyctor} is called a
\emph{factor block circulant matrix} where the factor is  $\diag(t^n,\ldots, t^n)$.
It is known that a factor block circulant matrix admits a block-diagonalization
in terms of the \emph{representer} $P(z)= W_0+W_1 z +\cdots + W_{n-1}z^{n-1}$:
\be
\begin{pmatrix*}[c]
	W_0 &  W_1& \cdots &  W_{n-2} & W_{n-1} \\
	t^{n} W_{n-1} &  W_0& \cdots&  W_{n-3} & W_{n-2}  \\
	\vdots & \vdots & \ddots & \vdots \\
	t^{n} W_1 & t^n W_2& \cdots & t^n  W_{n-1} &W_0 
\end{pmatrix*} = V \begin{pmatrix*}[c]
 P(t) &   &  &   & \\
 &   P(\omega \, t) &  &   &   \\
 &  & \ddots &  \\
  & &  &  & P(\omega^{n-1} \, t) 
\end{pmatrix*} V^{-1}
\ee
where $\omega$ is a primitive $n$-th root of unity and $V$ is a block Vandermonde
matrix
\be
V= \begin{pmatrix}
I & I &\cdots  & I \\
H_0  &H_1 &\cdots & H_{n-1} \\
H_0^2  &H_1^2 &\cdots & H_{n-1}^2 \\
\vdots & \vdots & & \vdots \\
H_0^{n-1}  &H_1^{n-1} &\cdots & H_{n-1}^{n-1}\\		
\end{pmatrix}, \quad H_i = w^{i} \, \diag(t,\ldots,t).
\ee
We refer to~\cite{circulant} for details. In particular, we have
\be	
\det \begin{pmatrix*}[c]
W_0 &  W_1& \cdots &  W_{n-2} & W_{n-1} \\
t^{n} W_{n-1} &  W_0& \cdots&  W_{n-3} & W_{n-2}  \\
\vdots & \vdots & \ddots & \vdots \\
t^{n} W_1 & t^n W_2& \cdots & t^n  W_{n-1} &W_0 
\end{pmatrix*} =\prod_{\omega^n=1} \det P(\omega\, t)
\ee
where, by definition,
\begin{align*}
P(\omega\, t) &= \sum_{r=0}^{n-1} W_r (\omega\, t)^r \\
&=\sum_{r=0}^{n-1} \sum_{k \in \BZ} \left( \sum_\square
  \mathbf{G}_{nk+r}^{\square} \mathrm{diag} (\zeta^\square) \right) t^{nk+r} \omega^{nk+r}
\\
&=\sum_{k \in \BZ} \left(\sum_\square \mathbf{G}_{k}^{\square}
  \mathrm{diag} (\zeta^\square)\right)(\omega\, t)^{k}
\qquad (\text{since} \,\, \omega^n=1).
\end{align*}
Combining the above calculations, we obtain
\be 
\tau^{\CS}(\calT^{(n)},t^n)  =
\prod_{\omega^n=1} \frac{\det\left(\sum_{k\in\BZ}
  \sum_\square  \mathbf{G}_{k}^{\square}\,
      \mathrm{diag} (\zeta^\square)\,(\omega\, t)^{k}\right)}{
  \prod_{j=1}^N\zeta_j^{f_j} \zeta_j'^{f_j'} \zeta_j''^{f_j''}} =
\prod_{\omega^n=1} \tau^{\CS}(\calT, \omega\, t).
\ee
\end{proof}

\begin{proof}(of Theorem~\ref{thm.2})
Let $R(t) = \sum_{\square} \mathbf{G}^\square (t)\, \mathrm{diag}(\zeta^\square)$ and
$R_i(t)$ denote its $i$-th row. Since the sum of all rows of  $\mathbf{G}^\square(1) =
\mathbf{G}^\square$ is $(2,\ldots,2)$ and
$\zeta_j + \zeta'_j+ \zeta''_j =0$, the sum $R_1(1)+\cdots+R_N(1)$ is
a zero vector. It follows that $\det R(1)=0$ and thus $\tau^{\CS}(\calT,1)=0$.	

Since the determinant is a linear map, we have
\be
\left. \frac{\partial}{\partial t} \right |_{t=1} \det R(t) 
= \sum_{i=1}^{N} \det
\begin{pmatrix}
R_1(1) \\
\vdots \\
\left.\frac{\partial}{\partial t}\right |_{t=1} R_i(t) \\
\vdots \\
R_N(1)
\end{pmatrix}  .
\ee	
For $i \neq N$ we exchange the $i$-th row
$\left.\frac{\partial}{\partial t} \right |_{t=1} R_i(t)$ with the $N$-th row
$R_N(1)$  and replace $R_N(1)$ by $-R_1(1)-\cdots-R_{N-1}(1)$. Then after some
elementary row operations, we obtain
\be
\left. \frac{\partial}{\partial t} \right |_{t=1} \det R(t) 
=\sum_{i=1}^N \det
\begin{pmatrix*}[c]
R_1(1) \\
\vdots \\
R_{N-1}(1) \\
\left.\frac{\partial}{\partial t}\right |_{t=1} R_i(t)
\end{pmatrix*}
= \det
\begin{pmatrix*}[c]
R_1(1) \\
\vdots \\
R_{N-1}(1) \\
 \left.\frac{\partial}{\partial t}\right |_{t=1} \sum_{i=1}^N R_i(t)
\end{pmatrix*}  .
\ee		
Comparing the definitions of equation~\eqref{1loop} and~\eqref{1loopt}, it
suffices to show that $\left.\frac{\partial}{\partial t}\right
|_{t=1} \sum_{i=1}^N R_i(t)$ agrees with
$\frac{1}{2}\sum_\square \mathbf{C}_\lambda^\square  \diag(\zeta^\square)$ up to
linear combination of $R_1(1),\ldots, R_{N}(1)$. 

Let $\ti \Sigma$ be the pre-image  of the peripheral torus $\Sigma$ under the
covering map $\ti M \rightarrow M$. Note that $\calT$ induces a triangulation
of $\ti \Sigma$ where vertices and triangles of $\ti \Sigma$ correspond to edges
and tetrahedra of $\ti \calT$, respectively. For each triangle $\Delta$ of
$\ti \Sigma$ corresponding to the tetrahedron $\mu^k \cdot \ti \Delta_j$, we
assign a \emph{corner parameter} $t^k \log z^\square_j$ to each corner of $\Delta$
according to the quad type of $\Delta_j$. Then the sum $r_i(t)$ of corner parameters
around a vertex of $\ti \Sigma$ corresponding to the edge $\ti e_i$ is
\be
\label{eqn.corner}
r_i(t)=\sum_{j=1}^N \sum_\square (\mathbf{G}^\square(t))_{ij}
\log z_j^\square .
\ee
It follows that (see Equation~\eqref{eqn.zeta})
\be
\label{eqn.Rr}
R_i(t) = \begin{pmatrix}
\frac{\partial}{\partial z_1} r_i(t) & \cdots &\frac{\partial}{\partial z_N} r_i(t)
\end{pmatrix}  .
\ee

We choose any lift $\ti \lambda$ of $\lambda$ to $\ti \Sigma$.
Since $\alpha(\lambda)=0$, $\ti \lambda$ is still a loop in $\ti \Sigma$.
Up to homotopy, we may assume that $\ti \lambda$ is an edge-path in $\ti \Sigma$
so that the surface $S$ bounded by $\ti \lambda$ and $\mu \cdot \ti \lambda$ is
triangulated. Here we view $S$ as a closed cylinder whose boundary is consisted of
$\ti \lambda$ and $\mu\cdot \ti \lambda$. Let $c_\lambda(t)$ (resp.,
$\check{c}_\lambda(t)$) be the sum of corner parameters in $S$ that are adjacent
(resp., not adjacent) to $\ti \lambda$. By definition, the equation $c_\lambda(1)=0$
represents the completeness equation of $\lambda$. Hence, 
\be
\label{eqn.lambda}
\sum_\square \mathbf{C}_\lambda^\square  \diag(\zeta^\square) = 
\begin{pmatrix}
  \frac{\partial}{\partial z_1} c_\lambda(1) & \cdots &
  \frac{\partial}{\partial z_N} c_\lambda(1)
\end{pmatrix}  .
\ee
Also, since $S$ consists of triangles and $\zeta_j + \zeta'_j +\zeta''_j=0$ for
all $j$, we have
\be
\label{eqn:ab}
\frac{\partial}{\partial z_j} (c_\lambda(t)+\check{c}_\lambda(t))=0  .
\ee
On the other hand, $t c_\lambda(t)$ is the sum of corner parameters  adjacent to
$\mu \cdot \ti \lambda$ and not in $S$. It follows that the sum of corner parameters
around vertices in $D \setminus \ti \lambda$ is given by
$tc_\lambda(t)+\check{c}_\lambda(t)$. Since $D \setminus \ti \lambda$ is a fundamental
domain of  $\Sigma$, exactly two vertices in $D \setminus \ti \lambda$  correspond
to a lift of $e_i$ for each $1 \leq i \leq N$. It follows that
\be 
t c_\lambda(t)+\check{c}_\lambda(t)= \sum_{i=1}^N (t^{a_i}+t^{b_i}) \, r_i(t)  .
\ee 	
for some $a_i$ and $b_i \in \BZ$ and thus
\begin{align}
  \left. \frac{\partial}{\partial t} \right|_{t=1}
  (t c_\lambda(t)+\check{c}_\lambda(t)) = \sum_{i=1}^N	(a_i+b_i) r_i(1) +2 \left.
    \frac{\partial}{\partial t} \right|_{t=1} \sum_{i=1}^Nr_i(t)  .
\end{align} 	
Taking the partial derivative with respect to $z_j$, we obtain
\begin{align}
  \left. \frac{\partial}{\partial t} \right|_{t=1} \frac{\partial}{\partial z_j}
  \sum_{i=1}^Nr_i(t) &= \frac{1}{2} \frac{\partial}{\partial z_j} \left.
    \frac{\partial}{\partial t} \right|_{t=1}
  (t c_\lambda(t)+\check{c}_\lambda(t))
  - \frac{1}{2} \sum_{i=1}^N (a_i+b_i) \frac{\partial}{\partial z_j}r_i(1)
  \\
  &=\frac{1}{2} \frac{\partial}{\partial z_j} c_\lambda(1)
  - \frac{1}{2} \sum_{i=1}^N (a_i+b_i) \frac{\partial}{\partial z_j}r_i(1) \,,
  \label{eqn.fin}
\end{align}
where the second equation is followed from the equation~\eqref{eqn:ab} with the
fact $t c_\lambda(t)+\check{c}_\lambda(t)=(t-1)c_\lambda(t)
+ c_\lambda(t)+\check{c}_\lambda(t)$. Combining the equations~\eqref{eqn.Rr},
\eqref{eqn.lambda}, and~\eqref{eqn.fin}, we obtain the desired result that
$\left.\frac{\partial}{\partial t}\right |_{t=1} \sum_{i=1}^N R_i(t)$ agrees with
$\frac{1}{2}\sum_\square \mathbf{C}_\lambda^\square  \diag(\zeta^\square)$ up to linear
combination of $R_1(1),\ldots, R_{N}(1)$:
\[
  \left. \frac{\partial}{\partial t} \right|_{t=1} \sum_{i=1}^N R_i(t)
  = \frac{1}{2}\sum_\square \mathbf{C}_\lambda^\square  \diag(\zeta^\square)
  + \frac{1}{2} \sum_{i=1}^N (a_i+b_i) R_i(1)  .
\]
\end{proof}

\begin{proof}(of Corollary~\ref{cor.sym})
The statement is equivalent to  
\be
\label{eqn.symm}
\det\left( \mathbf{A}(t)\, \mathrm{diag}(\zeta) + \mathbf{B}(t)\, \mathrm{diag}(\zeta'')
\right)= \ve t^r
\det\left( \mathbf{A}(1/t)\, \mathrm{diag}(\zeta) + \mathbf{B}(1/t)\,
\mathrm{diag}(\zeta'') \right) 
\ee
for some $\ve = \pm1$ and some integer $r$ and some choice of quads of $\calT$.
We make such a choice such that $\mathbf{B}$ is non-singular. This is always possible;
see~\cite[App.A]{DG1}. Then, $\det \mathbf{B}(t) \neq 0$, since
$\mathbf{B}(1)=\mathbf{B}$. Conjecture~\ref{conj.ABt} implies that
$\det \mathbf{B}(t)=\ve t^r \det \mathbf{B}(1/t) \neq 0$.
This, combined with Theorem~\ref{thm.sympl} gives
\begin{align*}
&\det\left( \mathbf{A}(t) \,\mathrm{diag}(\zeta)
+ \mathbf{B}(t)\, \mathrm{diag}(\zeta'') \right)\\
&= \det \mathbf{B}(t) \det \left(\mathbf{B}(t)^{-1} \mathbf{A}(t)
+ \mathrm{diag}(\zeta''/\zeta)\right) \det(\mathrm{diag}(\zeta))
\\
& =\ve t^r\det \mathbf{B}(1/t) \det \left(\mathbf{B}(1/t)^{-1} \mathbf{A}(1/t)
+ \mathrm{diag}(\zeta''/\zeta)\right) \det(\mathrm{diag}(\zeta))
\\
&=\ve t^r\det\left( \mathbf{A}(1/t)\, \mathrm{diag}(\zeta) + \mathbf{B}(1/t)\,
\mathrm{diag}(\zeta'') \right)
\end{align*}
Note that Theorem~\ref{thm.sympl} implies that
$\mathbf{B}(t)^{-1} \mathbf{A}(t)= (\mathbf{B}(1/t)^{-1} 
\mathbf{A}(1/t))^T$ which is used in the second equality above. 
\end{proof}

We end this section by giving a proof of Equation~\eqref{behavior2}. Taking the
derivative of Equation~\eqref{1loopcyclic} at $t=1$ and then simplifying
the result by Theorem~\ref{thm.2}, we obtain
\be
\label{behavior1}
\frac{\tau^\CS_\lambda(\calT^{(n)})}{\tau^\CS_\lambda(\calT)}
= \frac{1}{n}  \prod_{ \substack{\omega^n=1 \\ \omega \neq 1}} \tau^\CS(\calT,w)
\ee
which is equivalent to Equation~\eqref{behavior2} 
due to Theorem 4.1 of \cite{Porti:torsion} (explicitly, the ratio
$\tau_\lambda^{\CS}(\calT)/\tau_\mu^{\CS}(\calT)$ is the cusp shape of $M$ and
$\tau_\lambda^{\CS}(\calT^{(n)})/\tau_\mu^{\CS}(\calT^{(n)})$ is $n$ times of that).


\section{Examples, continued}
\label{sec.examplesb}

In this section, which is a continuation of Section~\ref{sec.examples},
we compute the twisted 1-loop invariant of an ideal triangulation
using the methods of \texttt{SnapPy}.

\subsection{The  knot $4_1$}
\label{sub.41}


In Section~\ref{sub.41a} we already gave the twisted NZ matrices of
the $4_1$ knot. Using the same notation, we now give the remaining data, namely
the shapes of the complete hyperbolic structure and flattenings, which are needed
to compute the twisted 1-loop invariant.

The solution for the complete structure is $z_1=z_2=\frac{1+\sqrt{-3}}{2}$.
Equations~\eqref{Gt41} and~\eqref{ABt41} give
\be
\mathbf{G}(t)\, \diag(\zeta) + 
\mathbf{G}'(t)\, \diag(\zeta') + 
\mathbf{G}''(t)\, \diag(\zeta'') = 
\begin{pmatrix}
2\zeta_1 t + \zeta_1' t^2 & \zeta_2' + 2 \zeta_2 t\\
\zeta_1' t+ 2 \zeta_1'' t^2 & 2 \zeta_2'' t + \zeta_2' t^2
\end{pmatrix}  .
\ee
One easily checks that $(f_1,f_2)=(0,0)$, $(f'_1,f'_2)=(1,1)$, and
$(f''_1,f''_2)=(0,0)$ are a combinatorial flattening of $\calT$.
Therefore, the twisted 1-loop invariant of $\calT$ is given by
\begin{align}
\tau^{\CS}(\calT,t)&= \pm \frac{\det
\begin{pmatrix}
2\zeta_1 t + \zeta_1' t^2 & \zeta_2' + 2 \zeta_2 t\\
\zeta_1' t+ 2 \zeta_1'' t^2 & 2 \zeta_2'' t + \zeta_2' t^2
\end{pmatrix}
}{\zeta'_1 \zeta'_2}\\
&= \pm t(t-1) \left(t^2- \frac{z_1 z_2+2z_1 +2z_2-4}{z_1z_2}\, t +1\right)  .
\end{align}
Substituting the solution for the complete structure into the above, we obtain
\[
  \tau^{\CS}(\calT,t)=\pm t(t-1)(t^2-5t+1)
\]
well-defined up to multiplication by a monomial $\pm t^r$ for some integer $r$.

\subsection{The  knot $6_3$}
\label{sub.63}

In this section we give the geometric solution and the flattenings of the triangulation
of the $6_3$ knot to compute its twisted 1-loop invariant, following the notation
of Section~\ref{sub.63a}.

The trace field of the $6_3$ knot is the number field of type $[0,3]$ and discriminant
$-11 \cdot 31^2$ given by $F=\BQ[\xi]$ where $\xi \approx 1.073 - 0.558 \sqrt{-1}$
satisfies
$\xi^6 - \xi^5 - \xi^4 + 2 \xi^3 - \xi + 1=0$.  
The solution for the complete structure is given exactly
\be
(z_1, z_2, z_3, z_4, z_5, z_6)
=(-\xi^2 + \xi, -\xi^2 + 1, -\xi^3 + \xi,-\xi^2 + 1, -\xi^3 + \xi,-\xi^2 + \xi)
\ee
and approximately by
\begin{center}
\begin{tabular}{lll}
  $z_1 \approx 0.23279 + 0.64139 \, i$,
  & $ z_2\approx 0.15884 + 1.20014 \, i $,
  & $z_3 \approx 0.84116 + 1.20014 \, i$, \\
  $z_4 \approx  0.15884 + 1.20014 \, i $,
  & $z_5  \approx 0.84116+ 1.20014 \, i$,
  &  $z_6 \approx 0.23279 + 0.64139 \, i$  .
\end{tabular}
\end{center}

It is easy to check that $(f_1,\ldots,f_6)=(0,1,0,1,0,0)$,
$(f'_1,\ldots,f'_6)=(1,0,1,0,1,1)$, and $(f''_1,\ldots,f''_6)=(0,0,0,0,0,0)$ is
a combintorial flattening of $\calT$. Equations~\eqref{Gt63} and~\eqref{ABt63}
give that
\begin{align*}
&\tau^{\CS}(\calT,t) =\frac{\det \left( \mathbf{G}(t)\, \diag(\zeta) + 
	\mathbf{G}'(t)\, \diag(\zeta') + 
	\mathbf{G}''(t)\, \diag(\zeta'')\right)
}{\zeta'_1 \zeta_2 \zeta_3' \zeta_4 \zeta'_5 \zeta'_6}\\
&= \frac{\det
\begin{pmatrix}
\zeta_1 t &  \zeta_2 t^2 & \zeta_3 t & \zeta_4 & 0 & \zeta_6 t\\
\zeta_1' & \zeta_2'' t + \zeta_2' t^2 + \zeta_2''t^3 & 0
& \zeta_4''+\zeta_4' t + \zeta_4''t^2  & 0 & \zeta_6' t^3  \\
\zeta_1'' & \zeta_2' t & \zeta_3't & 0 & \zeta_5' & \zeta_6'' t \\
\zeta_1'' t^2 & 0 & \zeta_3' t^2 & \zeta_4' t^2 & \zeta_5' t^3 & \zeta_6'' t^3 \\
\zeta_1' t^2 & 0 & \zeta_3'' + \zeta_3 t + \zeta_3''t^2 & 0
& \zeta_5'' +\zeta_5 t + \zeta_5'' t^2 & \zeta_6'\\
\zeta_1 t^2 & \zeta_2 t^2 & 0& \zeta_4 t^2 & \zeta_5 t^2 & \zeta_6  t^2
\end{pmatrix}
}{\zeta'_1 \zeta_2 \zeta_3' \zeta_4 \zeta'_5 \zeta'_6}\\
\end{align*}
Substituting the solution for the complete structure into the above, we obtain
\begin{align*}
  \tau^{\CS}(\calT,t) = & \, (t-1) t^8 (44 - 15 \xi - 15 \xi^2 + 34 \xi^3 - 19 \xi^5 \\
  & +(t+t^{-1})(-31 + 7 \xi + 7 \xi^2 - 10 \xi^3 + 3 \xi^5) \\
  & +(t^2+t^{-2})(15 - 3 \xi - 3 \xi^2 + 2 \xi^3 + \xi^5)
  -5(t^3+t^{-3}) + (t^4+t^{-4}))
\\
\approx & \, t^4 (-1.000 + 6.000 t - 12.805 t^2 + 33.472 t^3
- 85.242 t^4\\
 & + 85.242 t^5 - 33.472 t^6 + 12.805 t^7 -6.000 t^8 + 1.000 t^9)  .
\end{align*}
up to multiplication by a monomial $\pm t^r$ for some integer $r$.

\subsection*{Acknowledgments}
The first author wishes to thank Nathan Dunfield for enlightening conversations
and for sharing his drawing program of the triangulated cusp of an ideal
triangulation. The second author was supported by Basic Science Research Program
through the NRF of Korea funded by the Ministry of Education (2020R1A6A3A03037901).


\bibliographystyle{hamsalpha}
\bibliography{biblio}

\providecommand{\bysame}{\leavevmode\hbox to3em{\hrulefill}\thinspace}
\providecommand{\href}[2]{#2}
\providecommand{\eprint}{\begingroup \urlstyle{rm}\Url}
\begin{thebibliography}{RCdSL88}

\bibitem[AGK]{AGK:KLV}
J{\o}rgen~Ellegaard Andersen, Stavros Garoufalidis, and Rinat Kashaev,
  \emph{The volume conjecture for the klv state-integral}, Preprint 2021.

\bibitem[CDW]{snappy}
Marc Culler, Nathan Dunfield, and Jeffrey Weeks, \emph{Snap{P}y, a computer
  program for studying the topology of $3$-manifolds}, Available at
  \url{http://snappy.computop.org} (30/01/2015).

\bibitem[Cho06]{choi06}
Young-Eun Choi, \emph{Neumann and {Z}agier's symplectic relations}, Expo. Math.
  \textbf{24} (2006), no.~1, 39--51.

\bibitem[Dav79]{Davis}
Philip Davis, \emph{Circulant matrices}, A Wiley-Interscience Publication, John
  Wiley \& Sons, New York-Chichester-Brisbane, 1979.

\bibitem[DFJ12]{Dunfield:twisted}
Nathan Dunfield, Stefan Friedl, and Nicholas Jackson, \emph{Twisted {A}lexander
  polynomials of hyperbolic knots}, Exp. Math. \textbf{21} (2012), no.~4,
  329--352.

\bibitem[DG13]{DG1}
Tudor Dimofte and Stavros Garoufalidis, \emph{The quantum content of the gluing
  equations}, Geom. Topol. \textbf{17} (2013), no.~3, 1253--1315.

\bibitem[DG16]{Dubois-Garoufalidis}
Jerome Dubois and Stavros Garoufalidis, \emph{Rationality of the {${\rm
  SL}(2,\Bbb C)$}-{R}eidemeister torsion in dimension 3}, Topology Proc.
  \textbf{47} (2016), 115--134.

\bibitem[DG18]{DG2}
Tudor Dimofte and Stavros Garoufalidis, \emph{Quantum modularity and complex
  {C}hern-{S}imons theory}, Commun. Number Theory Phys. \textbf{12} (2018),
  no.~1, 1--52.

\bibitem[DGG13]{DGG2}
Tudor Dimofte, Davide Gaiotto, and Sergei Gukov, \emph{3-manifolds and 3d
  indices}, Adv. Theor. Math. Phys. \textbf{17} (2013), no.~5, 975--1076.

\bibitem[DGG14]{DGG1}
\bysame, \emph{Gauge theories labelled by three-manifolds}, Comm. Math. Phys.
  \textbf{325} (2014), no.~2, 367--419.

\bibitem[DY12]{DubYam}
J\'{e}r\^{o}me Dubois and Yoshikazu Yamaguchi, \emph{The twisted {A}lexander
  polynomial for finite abelian covers over three manifolds with boundary},
  Algebr. Geom. Topol. \textbf{12} (2012), no.~2, 791--804.

\bibitem[Fox56]{Fox3}
Ralph Fox, \emph{Free differential calculus. {III}. {S}ubgroups}, Ann. of Math.
  (2) \textbf{64} (1956), 407--419.

\bibitem[GK19]{GK:meromorphic}
Stavros Garoufalidis and Rinat Kashaev, \emph{A meromorphic extension of the
  3{D} index}, Res. Math. Sci. \textbf{6} (2019), no.~1, Paper No. 8, 34.

\bibitem[Kit96]{Kitano96}
Teruaki Kitano, \emph{Twisted {A}lexander polynomial and {R}eidemeister
  torsion}, Pacific J. Math. \textbf{174} (1996), no.~2, 431--442.

\bibitem[KL99]{KL99}
Paul Kirk and Charles Livingston, \emph{Twisted {A}lexander invariants,
  {R}eidemeister torsion, and {C}asson-{G}ordon invariants}, Topology
  \textbf{38} (1999), no.~3, 635--661.

\bibitem[KLV16]{KLV}
Rinat Kashaev, Feng Luo, and Grigory Vartanov, \emph{A {TQFT} of
  {T}uraev-{V}iro type on shaped triangulations}, Ann. Henri Poincar\'e
  \textbf{17} (2016), no.~5, 1109--1143.

\bibitem[Neu92]{Neumann04}
Walter Neumann, \emph{Combinatorics of triangulations and the {C}hern-{S}imons
  invariant for hyperbolic {$3$}-manifolds}, Topology '90 ({C}olumbus, {OH},
  1990), Ohio State Univ. Math. Res. Inst. Publ., vol.~1, de Gruyter, Berlin,
  1992, pp.~243--271.

\bibitem[NZ85]{NZ}
Walter Neumann and Don Zagier, \emph{Volumes of hyperbolic three-manifolds},
  Topology \textbf{24} (1985), no.~3, 307--332.

\bibitem[Por97]{Porti:torsion}
Joan Porti, \emph{Torsion de {R}eidemeister pour les vari\'{e}t\'{e}s
  hyperboliques}, Mem. Amer. Math. Soc. \textbf{128} (1997), no.~612, x+139.

\bibitem[RCdSL88]{circulant}
Julio~Cesar Ruiz-Claeyssen and Liara~Aparecida dos Santos~Leal,
  \emph{Diagonalization and spectral decomposition of factor block circulant
  matrices}, Linear Algebra Appl. \textbf{99} (1988), 41--61.

\bibitem[Sie]{siejakowski}
Rafa{\l} Siejakowski, \emph{Infinitesimal gluing equations and the adjoint
  hyperbolic {R}eidemeister torsion}, \eprint{arXiv:1710.02109}, Preprint 2017.

\bibitem[Thu77]{Thurston}
William Thurston, \emph{The geometry and topology of 3-manifolds},
  Universitext, Springer-Verlag, Berlin, 1977,
  \url{http://msri.org/publications/books/gt3m}.

\bibitem[Wad94]{Wad94}
Masaaki Wada, \emph{Twisted {A}lexander polynomial for finitely presentable
  groups}, Topology \textbf{33} (1994), no.~2, 241--256.

\bibitem[Yam08]{Yamaguchi}
Yoshikazu Yamaguchi, \emph{A relationship between the non-acyclic
  {R}eidemeister torsion and a zero of the acyclic {R}eidemeister torsion},
  Ann. Inst. Fourier (Grenoble) \textbf{58} (2008), no.~1, 337--362.

\bibitem[Zic16]{Zickert:Apoly}
Christian Zickert, \emph{Ptolemy coordinates, {D}ehn invariant and the
  {$A$}-polynomial}, Math. Z. \textbf{283} (2016), no.~1-2, 515--537.

\end{thebibliography}
\end{document}